\documentclass[10pt,english]{amsart}
\pdfoutput=1
\usepackage[T1]{fontenc}
\usepackage[latin9]{inputenc}
\usepackage{geometry}
\usepackage{babel}
\usepackage{refstyle}
\usepackage{amsmath}
\usepackage{amsthm}
\usepackage{amsbsy}
\usepackage{amstext}
\usepackage{graphicx}
\usepackage{esint}
\usepackage[normalem]{ulem}
\usepackage{color}
\usepackage{cancel}
\usepackage[unicode=true,pdfusetitle,
 bookmarks=true,bookmarksnumbered=false,bookmarksopen=false,
 breaklinks=false,pdfborder={0 0 1},backref=false,colorlinks=false]
 {hyperref}

\makeatletter

\AtBeginDocument{\providecommand\secref[1]{\ref{sec:#1}}}
\AtBeginDocument{\providecommand\eqref[1]{\ref{eq:#1}}}
\AtBeginDocument{\providecommand\figref[1]{\ref{fig:#1}}}

\newcommand{\lyxdot}{.}
  \definecolor{dkgr}{rgb}{0.0,0.6,0.2}
  \newcommand{\rrr}[1]{}
  \renewcommand{\rrr}[1]{{\bf\color{dkgr}#1}} 
  \newcommand{\?}{\hspace*{0.1em}}
  \newtheorem{note}{\textbf{Note}}
  \newcommand{\TheoremEnd}{~\hfill$\clubsuit\/$}
  \newcommand{\M}{\mathcal{M}}
  \newcommand{\R}{\mathcal{R}}
  \newcommand{\calL}{\mathcal{L}}
\RS@ifundefined{subref}
  {\def\RSsubtxt{section~}\newref{sub}{name = \RSsubtxt}}
  {}
\RS@ifundefined{thmref}
  {\def\RSthmtxt{theorem~}\newref{thm}{name = \RSthmtxt}}
  {}
\RS@ifundefined{lemref}
  {\def\RSlemtxt{lemma~}\newref{lem}{name = \RSlemtxt}}
  {}

\numberwithin{equation}{section}
\numberwithin{figure}{section}

\@ifundefined{showcaptionsetup}{}{%
 \PassOptionsToPackage{caption=false}{subfig}}
\usepackage{subfig}
\makeatother

\begin{document}

\title{Study of a model equation in detonation theory: multidimensional effects}
\thanks{This work was presented in part at the 25th International Colloquium on the Dynamics of Explosions and Reactive Systems (ICDERS) held at the University of Leeds, UK, on Aug. 2-7, 2015} 

\author{L. M. Faria$^{1,2}$, A. R. Kasimov$^{1}$,}
\thanks{$^{1}$Applied Mathematics and Computational Sciences, KAUST, Saudi Arabia ({\tt aslan.kasimov@kaust.edu.sa}). L.F.'s present address: Department of Mathematics, MIT, USA  ({\tt lfaria@mit.edu}). L. F. and A. K. gratefully acknowledge research support by King Abdullah University of Science and Technology (KAUST)}
\author{R. R. Rosales$^{2}$}
\thanks{$^{2}$Department of Mathematics, MIT, USA  ({\tt rrr@math.mit.edu}). The work by R.R.R. was partially supported by NSF grants DMS-1115278 and DMS-1318942.}

\maketitle
\begin{abstract}
We extend the reactive Burgers equation presented in
\cite{Kasimov2013,Faria2014} to include multidimensional effects. Furthermore, we explain
how the model can be rationally justified
following the ideas of the
asymptotic theory developed in \cite{Faria2014b}.
The proposed model is a forced version of the unsteady small disturbance
transonic flow equations. We show that for physically reasonable choices of
forcing functions, traveling wave solutions akin to detonation waves exist.
It is demonstrated that multidimensional effects play an important role
in the stability and dynamics of the traveling waves. Numerical simulations indicate
that solutions of the model tend to form multi-dimensional patterns
analogous to cells in gaseous detonations.
\end{abstract}
%


\pagestyle{myheadings}
\thispagestyle{plain}
\markboth{L. M. Faria, A. R. Kasimov, and R. R. Rosales}{Study of a model equation in detonation theory: multidimensional effects}

\maketitle

\section{Introduction}
The fact that very simple mathematical models, such as the logistic
map, can produce extremely complicated solutions was a surprising
and counterintuitive discovery in mathematics \cite{may1976simple}.
It suggested that complex behavior need not always be described by
complicated equations. The hope with such simple models is that they
can shed light into fundamental mechanisms responsible for complexity,
while discarding secondary details which only obscure the important
underlying dynamics. In this work, we introduce a simple system of
partial differential equations, consisting of a non-locally forced
version of the unsteady transonic small disturbance equations
\cite{lin1948two} (also closely related to Kadomtsev-Petviashvili
equation in water waves \cite{kadomtsev1970stability} and Zabolotskaya-Khokhlov
equation in nonlinear acoustics \cite{Zab-Khokhlov-1969}). The system
is shown to reproduce some of the structure and dynamics of two-dimensional
gaseous detonations. 

Detonations are a type of combustion process in which strong pressure
waves ignite, and are sustained by, exothermic chemical reactions.
The equations of combustion theory pose a formidable challenge from
a theoretical point of view because they couple a compressible flow
description (Euler/Navier-Stokes equations) to chemical kinetics.
Indeed, the main difficulty in understanding most reactive flows lies
precisely in this two-way coupling: the wave initiates chemical reactions,
and the reactions sustain the wave, with either ceasing to exist without
the other. To make matters even more complicated, detonations tend
to be multi-dimensional and unsteady. It is thus not surprising that
models simpler than the reactive Euler/Navier-Stokes equations are
highly desirable for their understanding. 

The first attempt at a reduced qualitative description of detonations goes back to the
 work by Fickett \cite{Fickett1979,Fickett1985b}, who introduced a toy model (called
 an ``analog'' by Fickett) as a vehicle to yield a better understanding of the
 intricacies of detonation theory. Others took a similar 
qualitative approach. Majda, for example, focused on the effect of viscosity on
combustion waves, and showed through a simplified model that a theory analogous to
ZND theory (the classical inviscid theory of one-dimensional steady detonations due
to Zel'dovich \cite{Zeldovich1940}, von Neumann \cite{vonNeumann1942}, and D\"oring
\cite{Doering1943}) exists for viscous detonations \cite{Majda1981}. Radulescu and
Tang \cite{Radulescu2011} recently demonstrated that simplified models can capture
not only the steady states, but also much of the unsteady dynamics of one-dimensional
detonations. The importance of this latter development stems from the fact that
unsteadiness is a common feature of gaseous detonations. Along the same lines of
analog modeling, 
in \cite{Kasimov2013,Faria2014} we have shown
that even a scalar forced Burgers equation contains all the ingredients necessary to
reproduce the complexity of one-dimensional detonations, including complicated chaotic
solutions. 

Importantly, all the prior work with analog models
is limited to one-dimensional descriptions. Therefore, it cannot capture the important
effects played by transverse shock waves, which lead to the generation
of cellular patterns in gaseous detonations \cite{Fickett2012,Lee2008}.
The purpose of this paper is to extend our existing one-dimensional
model \cite{Kasimov2013,Faria2014} by introducing a multi-dimensional
analog of detonations. With this extended model we show, through
stability analysis and numerical simulations, that some multi-dimensional
detonation properties are amenable to descriptions significantly simpler than
the reactive Euler equations.

The remainder of this paper is organized as follows. In \secref{The-model-2d},
 the two-dimensional analog 
model is introduced, and some motivation based on the theory of weakly
curved hyperbolic waves is given. We also discuss 
how nonlocal forcing functions,
such as the one studied in \cite{Faria2014}, can arise through an
asymptotic approximation of a local reaction rate.
In \secref{Traveling-wave-solutions-2d}, we present traveling wave
solutions of the model, together with a solution of the linear stability problem by means
of the Laplace transform. Finally, in \secref{An-example-2d-toy} we illustrate
the linear stability theory and nonlinear dynamics by studying a concrete
example. We demonstrate by this particular example that: (1) transverse
perturbations are typically more unstable than purely longitudinal
ones, and (2) when unstable, the traveling wave solutions tend to
form multi-dimensional patterns of varying complexity, depending on
the distance from the neutral stability boundary.

\section{The two-dimensional analog}\label{sec:The-model-2d}
As a starting point, we take the one-dimensional model introduced
in \cite{Kasimov2013},
\begin{equation}\label{eq:toy-model-1d}
 u_{t} + \frac{1}{2}\left(u^{2}\right)_{x} = f\left(x-x_{s},u_{s}\right),
\end{equation}
where $u$ is the main state variable (playing the role of, say, the flow velocity),
 $x_{s}=x_s(t)$ denotes the position of the leading edge of the combustion front,
 $u_{s}=u_s(t)$ denotes the state immediately behind $x_s$, 
 $f$ is the reaction term, and the subscripts $t$ and $x$ denote time and space
 derivatives, respectively. There is some flexibility in the precise definition of
 $x_{s}$, especially when viscous effects (not considered here) are included. In this
 paper, we shall focus on the case where the leading edge of the combustion front is
 given by a discontinuity (shock) in $u$, ahead of which no chemical reaction takes
 place. Then, $x_{s}$ denotes the shock position, $u_{s}$ is the post-shock state, and
 \vspace*{-0.5em}
 \begin{equation}\label{eq:fnonnegint}
    f = 0\;\;\mbox{for}\;\;x>x_s.
 \end{equation}%
We also assume that $f$
    is non-negative and integrable, with $\int f\/d\/x > 0$.
    
In order to extend (\ref{eq:toy-model-1d}) to two spatial dimensions, we introduce
another variable, $v$, to describe the transverse velocity. Then, a relation between
$u$ and $v$ is required to close the system. As this is a qualitative theory,
we are free to extend the model as we like. Some ways, however, make more physical
sense than others. Our extension here is motivated by the asymptotic
 equations governing
weakly curved hyperbolic waves, in which the dependence on the transverse direction
is linear, and reinforces the fact that weakly nonlinear quasi-planar waves generate
no vorticity to leading order. We thus propose the following two-dimensional
extension of \eqref{toy-model-1d}:
\begin{eqnarray}
 u_{t} + \frac{1}{2}\left(u^{2}\right)_{x} + v_{y} & = & f(x-x_{s},u_{s}),
 \label{eq:toy-model-2d-1}\\
 v_{x} - u_{y} & = & 0, \label{eq:toy-model-2d-2}
\end{eqnarray}
 where now $x_s = x_s(y\/,\,t)$ and $u_s = u_s(y\/,\,t)$. By eliminating $v$, this system can also be written as a single equation, either of a second order, 
\begin{equation}
u_{xt} + \frac{1}{2}(u^{2})_{xx} + u_{yy}  =  f(x-x_{s},u_{s})_x, 
\end{equation}
 or as an integro-partial differential equation. However, we keep it as a system for the purpose of the analysis and numerical computations below.

In the absence of the source function $f$, equations
(\ref{eq:toy-model-2d-1}-\ref{eq:toy-model-2d-2})
 provide the
canonical description for weakly-nonlinear quasi-planar hyperbolic waves, and have
been derived in many different contexts in the past. For example, in the study of
flow past a body in a compressible fluid \cite{lin1948two}, and in nonlinear
acoustics \cite{Zab-Khokhlov-1969}. Furthermore, similar equations have also been
obtained in the study of water waves \cite{kadomtsev1970stability}, where a
dispersive term proportional to $u_{xxx}$ is included.
 An example in combustion theory can be found in \cite{rosales1989diffraction}.

 The (\ref{eq:toy-model-2d-1}-\ref{eq:toy-model-2d-2}) system is a natural extension
 of \eqref{toy-model-1d}. The new equations include both canonical weakly-nonlinear
 2D hyperbolic effects, as well as a source term which aims at capturing the effects
 of the energy released by chemical reactions.
 Additional dispersive or dissipative effects, even though quite interesting, are not
 studied here. This new system must be
interpreted as an \emph{ad-hoc} model. However,
 it is closely related to the asymptotic model in \cite{Faria2014thesis,Faria2014b},
 which can be obtained by a systematic reduction of the
 reactive Navier-Stokes equations. In order to see the connection, we recall some of
 the details presented in \cite{Faria2014thesis}. When dissipative effects (i.e. viscosity,
 heat conduction, and diffusion) are ignored, the asymptotic model 
 derived in \S 5 of
 \cite{Faria2014thesis} reduces to
\begin{eqnarray}
 u_{\tau} + u\?u_{x} + v_{y} & = & -\frac{1}{2}\?\lambda_{x},
    \label{eq:asymptotic-model-1}\\
 v_{x}  -  u_{y} & = &0,\label{eq:asymptotic-model-2}\\
 \lambda_{x} & = &
    \begin{cases}
      0 & \text{for }T<T_{i}\\
     -k\?(1-\lambda)\?\exp\left(\theta\?q\?T\right) & \text{for }T\geq T_{i}, 
    \end{cases}
    \label{eq:asymptotic-model-3}\\
 T & = & \frac{u}{\sqrt{q}}+\lambda,\label{eq:asymptotic-model-4}
\end{eqnarray}
where the dependent variables $u,\, v,\, T$, and $\lambda$ represent the leading
order perturbations to the $x$ velocity, $y$ velocity, temperature, and reaction
progress variable, respectively. The independent variables $x$ and $y$ represent
longitudinal and transverse spatial coordinates relative to a moving acoustic frame,
and $\tau$ is a slow time variable. The parameters $q,\,\theta,$ and $T_{i}$ are the
rescaled heat release, activation energy, and ignition temperature, respectively.
The heat release is a measure of how much chemical energy is contained in the
mixture, while the activation energy controls the sensitivity of the chemical
reactions to temperature. The ignition temperature, $T_{i}$, sets a threshold below which no chemical
 reactions take place --- i.e., the reaction rate $\omega$, defined by the right
 hand side in (\ref{eq:asymptotic-model-3}), vanishes.

In order to bridge (\ref{eq:toy-model-2d-1}-\ref{eq:toy-model-2d-2}) and the
asymptotic model (\ref{eq:asymptotic-model-1}-\ref{eq:asymptotic-model-4}),
 we must justify/motivate the replacement of the term $-\frac{1}{2}\,\lambda_x$ in
 (\ref{eq:asymptotic-model-1}), by the nonlocal forcing $f(x-x_{s},u_{s})$ present
 in (\ref{eq:toy-model-2d-1}).
 In previous work \cite{Kasimov2013} we justified this by  
 adding the extra assumption that 
 $\omega=\omega\left(\lambda,u_{s}\right)$ --- that is, we considered a reaction rate that
 depends solely on how strong the precursor shock wave is.
This simplifying assumption has also been used in the past, both in the context of
analog modeling \cite{Fickett1979,Fickett1985b} and condensed explosives, but without a rational
justification. We explain next how the nonlocal approximation can be justified as
a uniformly valid asymptotic approximation to $\omega$, when the scaled heat release $q$ is
sufficiently large. 

The key observation is that, when $q\gg1,$ the main contribution of $u$ in
(\ref{eq:asymptotic-model-4}) comes from the region where $\lambda\ll1$. In
particular, for $u$ to contribute to the reaction rate, we need
$\sqrt{q}\lambda\approx u$, which occurs when $\lambda=O(1/\sqrt{q})$. This means
that the reaction rate is appreciably affected by $u$ only when $\lambda\approx0$,
and since $\lambda=0$ ahead of the combustion front, we know that $\lambda$ is
small only near $x_{s}$, where $u$ is well approximated by $u_{s}$
(assuming that $u$ is sufficiently smooth to the left of $x_{s}$).
 
 More formally, let $\epsilon=1/\sqrt{q}$, and write 
 $T=\lambda+\epsilon u(\lambda\/, t)$. This is possible because, at any fixed time
    $t$, we can replace $x$ by $\lambda$ as a parameter, since $\lambda(x)$ is
    a monotone function of $x$ for $x\leq x_{s}$. From here, it would seem that a good
 approximation to the temperature is given by $T \sim T^{\text{nonunif}} = \lambda$, since
 $T = \lambda + O(\epsilon)$. However, this approximation breaks down for
 $\lambda \leq O(\epsilon)$, when the relative error becomes $O(1)$ or
 larger (note that $u\,\rule[-1.0ex]{0.3mm}{2.2ex}_{\lambda=0}=u_s(t)>0\/$
    because of the entropy condition satisfied by the lead shock at $x=x_s$).
 It is easy to see that an approximation uniformly valid for $\lambda\in[0,1]$ (in
 the sense that the relative error is small for all values of $\lambda$) is given by
 \begin{equation}\label{eq:Tunif}
   T \sim T^{\text{unif}} = \lambda+\epsilon\,u\,\rule[-1.0ex]{0.3mm}{2.2ex}_{\lambda=0}
   = \lambda + \epsilon\,u_s\/.
 \end{equation}
 This holds provided that $u$ as a function of $\lambda$ is smooth enough from the left near
 $\lambda = 0$, and bounded. Thus the approximation has an $O(\epsilon\,\lambda)$
 error, which yields an $O(\epsilon)$ relative error --- recall that $0 \leq \lambda$,
 $0 < \epsilon \ll 1$, and $u_s = O(1) > 0$. It is important to note that this
 approximation does not require $u$ to be well approximated by $u_s$ everywhere. For
 example, for the traveling wave profiles that we study later in
 \secref{Traveling-wave-solutions-2d}, $\max_{x<x_{s}}\left| u-u_s \right|=O(1)$, yet
 (\ref{eq:Tunif}) remains valid. 

 Substituting (\ref{eq:Tunif}) into (\ref{eq:asymptotic-model-3}) yields an equation
 of the form
 \[ \lambda_{x}=\omega\left(u_{s},\lambda\right), \quad \mbox{with}\;\,\lambda=0\,\;
    \mbox{at}\;\, x=x_s. \]
 This is an ODE for $\lambda$, and its solution has the form
 $\lambda=F(x-x_s\/,\,u_s)$. 
 Letting $$f(x-x_s,u_s) = - \frac{F_x(x-x_s, u_s)}{2},$$ and substituting this result into
 (\ref{eq:asymptotic-model-1}), 
 the model in (\ref{eq:toy-model-2d-1}-\ref{eq:toy-model-2d-2}) follows.

\begin{figure}[htb!]
 \begin{centering}
  \subfloat[$q=5$]{\includegraphics[width=3in,height=2in]{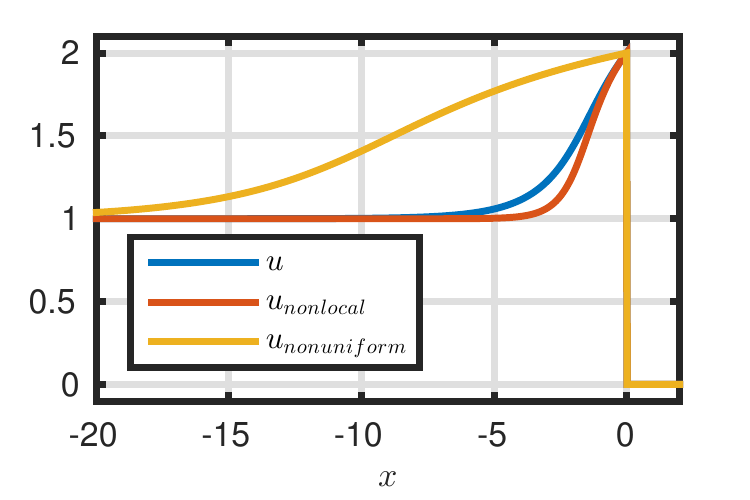}}
  \subfloat[$q=10$]{\includegraphics[width=3in,height=2in]{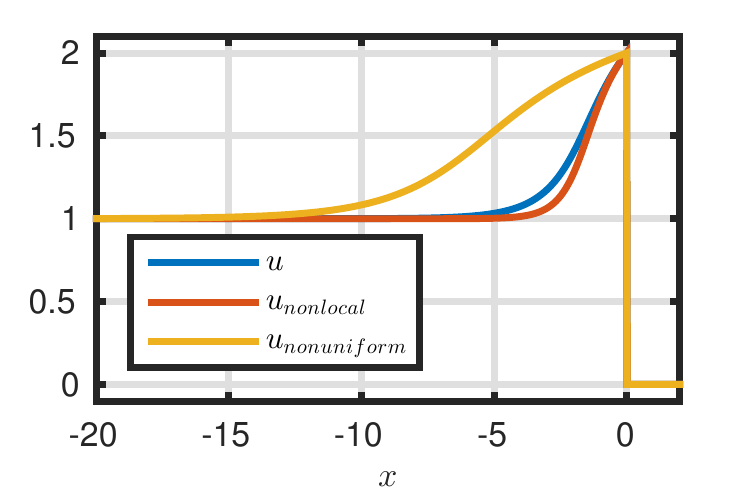}}
 \par\end{centering}

 \begin{centering}
  \subfloat[$q=100$]{\includegraphics[width=3in,height=2in]{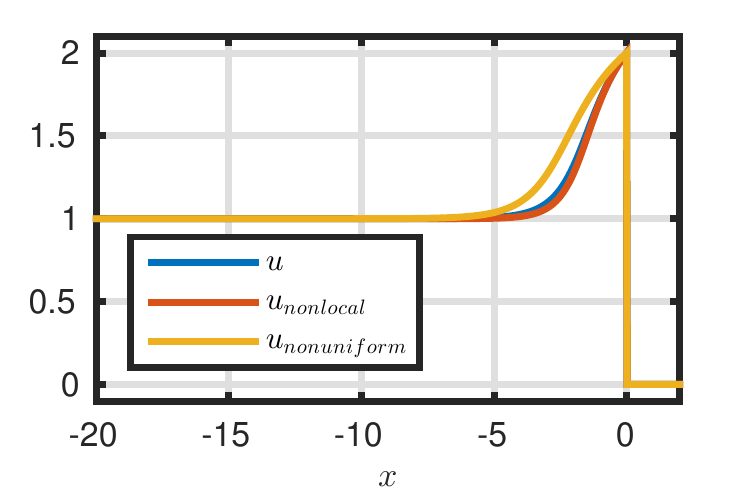}}
  \protect
  \caption{\label{fig:Local-vs.-nonlocal}%
     Comparisons of the steady-state solutions for
     (\ref{eq:asymptotic-model-1}-\ref{eq:asymptotic-model-4}), using
     the temperatures:
     (i)   $T=u/\sqrt{q}+\lambda$ (exact and local);
     (ii)  $T=u_{s}/\sqrt{q}+\lambda$ (uniform and nonlocal approximation); and
     (iii) $T=\lambda$ (nonuniform approximation).
  }\par\end{centering}
\end{figure}

 To illustrate the importance of uniformity in (\ref{eq:Tunif}), in
 \figref{Local-vs.-nonlocal}, we plot the one dimensional traveling wave solutions
 for the system in (\ref{eq:asymptotic-model-1}-\ref{eq:asymptotic-model-4}) using
 three different temperatures: (i) exact, given by the ``local'' formula
 $T = \lambda + \epsilon\,u$ in (\ref{eq:asymptotic-model-4}); (ii) uniform, given
 by the ``non-local'' formula in (\ref{eq:Tunif}); and (iii) non-uniform, given by
 $T = T^{\text{nonunif}} = \lambda$. The plots are for increasing values of $q$, with
 $q\,\theta=2$ kept fixed. Even for $q=5$ ($\epsilon=1/\sqrt{5}$) a reasonable
 agreement between the local and nonlocal profiles occurs. The agreement improves,
 as expected, for larger values of $q$. On the other hand, the nonuniform
 approximation has a significant departure from the local profile even for $q=100$.
%

One difficulty with the procedure outlined in the previous paragraph for
obtaining $f(x-x_{s},u_{s})$ from the reaction rate is that it may not be possible to
integrate $\lambda_{x}=\omega\left(u_{s}\/,\,\lambda\right)$ analytically, and in a sufficiently simple form. In the
physical example of an Arrhenius law, for instance,
 the equation to be integrated is
\begin{eqnarray}
 \lambda_{x} & = & -k\?(1-\lambda)\?\exp\left(\theta q\left(
    \frac{u_{s}}{\sqrt{q}} + \lambda\right)\right). \label{eq:nonlocal-rection-rate}
\end{eqnarray}
This yields $\lambda$ implicitly, in terms of the inverse of an exponential integral. This inverse
then needs to be differentiated in order to obtain $f(x-x_{s},u_{s})$, which produces a rather cumbersome expression. Because of these
complications, when focusing on a specific example later in \secref{An-example-2d-toy},
we opt for simplicity and use an \emph{ad hoc} explicit forcing function $f$ which
captures the main qualitative features of simple chemical reactions.
\section{Traveling wave solutions and stability analysis}
\label{sec:Traveling-wave-solutions-2d}
In this section, we investigate the one-dimensional traveling wave solutions of
(\ref{eq:toy-model-2d-1}-\ref{eq:toy-model-2d-2}),
 as well as 
their stability properties. Note that, since in the one-dimensional case the proposed
system reduces to \eqref{toy-model-1d}, the traveling wave solutions are identical
to those found in \cite{Faria2014}. For completeness, we briefly repeat the discussion
here.
\subsection{Traveling wave solutions}\label{sub:Traveling-wave-solutions-2d-toy}
We seek for solutions of (\ref{eq:toy-model-2d-1}-\ref{eq:toy-model-2d-2}) of
 the form $u(x,y,t)=u_{0}(x-Dt)$ and $v=0$ (more generally, $v$ constant) for
 $x\leq D\,t$, with $x_s=D\,t$ and $u_s=u_{s0}=u_0(0)$. It is easy to see that $u_0$
 must have the form
\begin{eqnarray}\label{eq:steady-state-2d-1}
 u_{0}(\xi) & = & D \pm \sqrt{C_* - 2\int_{\xi}^{0}f\left(z,u_{0s}\right)\?dz},
\end{eqnarray}
 where $\xi=x-D\?t$, and $C_*$ is an integration constant. The shock conditions
 at $x=x_s$, and the fact that the solution vanishes for $x > x_s$, yield
 $u_{s0} = 2D > 0$. This determines $C_*$ and the sign of the square root in
 (\ref{eq:steady-state-2d-1}). The full profile is then
\begin{equation}\label{eq:TWS-1d}
 v \equiv 0 \quad\mbox{and}\quad
 u(\xi) = \begin{cases}
   D+\sqrt{D^{2}-2\int_{\xi}^{0}f\left(z,u_{0s}\right)dz} & \text{for }\xi\leq0,\\
   0                                                   & \text{for }\xi>0.
 \end{cases}
\end{equation}
 These are the traveling waves for the model; they consist of a shock of strength
 $2D>0$ moving into an unperturbed, unburnt state. For these solutions to be real, it
 must be that
\begin{equation}\label{eq:u0s}
 D = \zeta\?\sqrt{2\int_{-\infty}^{0}f\left(y,u_{0s}\right)\?dy},
 \quad\mbox{where}\;\;\zeta \geq 1.
\end{equation}
The parameter $\zeta$ is the overdrive factor. The special value $\zeta=1$ corresponds
to the slowest possible wave speed -- the Chapman-Jouguet velocity. When $\zeta=1$, the
 characteristic speed (i.e., $u$) at the end of the reaction zone is equal to the
 wave speed $D$. On the other hand, waves with $\zeta>1$ are ``overdriven''.
 For them, the characteristic speed is everywhere greater than the wave speed. Because
 of this, overdriven detonations are not self-sustained: the characteristics from any
 $\zeta < 0$ catch up to the lead shock in a finite time, and affect its dynamics.
Furthermore, the larger the overdrive, the closer the profile resembles that of an
inert piston-induced shock. Therefore, in the limit of large overdrive, detonations
are expected to be stable. We also point out here that, as in \cite{Faria2014}, the total energy release is assumed to be constant, i.e., $\int f\/dx = q/2=const$. As a result, $D$ is given explicitly by (\ref{eq:u0s}). Without this constraint on $f$, (\ref{eq:u0s}) would have been an implicit equation for $D$, because $u_{0s} = 2D$.

Next we consider the multidimensional stability of the traveling waves given by
(\ref{eq:TWS-1d}).  
\subsection{Linear stability}\label{sub:Linear-stability-2d-toy}
We consider in this subsection the multi-dimensional linear stability for the
solutions given by (\ref{eq:TWS-1d}). For the purpose of the calculations that
follow, it is convenient to rewrite (\ref{eq:toy-model-2d-1}-\ref{eq:toy-model-2d-2})
in a shock-attached frame.
 Assume that the shock position is given by a single valued function,
 $x_s = s(y\/,\,t)$, and introduce the shock-attached variable $\xi=x-s(y,t)$ to
 replace $x$. Then (\ref{eq:toy-model-2d-1}-\ref{eq:toy-model-2d-2}) takes the form
\begin{eqnarray}
 u_{t} + (u-s_{t})\?u_{\xi} + v_{y} - s_{y}\?v_{\xi} & = & f(\xi,u_{s}),
    \label{eq:toy-model-2d-saf-1}\\
 u_{y} - s_{y}\?u_{\xi} - v_{\xi} & = & 0.\label{eq:toy-model-2d-saf-2}
\end{eqnarray}
The quantities $s_{t}$ and $s_{y}$ are related to the states at the shock by the jump
conditions, given by
\begin{eqnarray}
 s_{t}\left[u\right]-\frac{1}{2}\left[u^{2}\right]+s_{y}\left[v\right] & = & 0,
    \label{eq:rk1}\\
 s_{y}\left[u\right]+\left[v\right] & = & 0,\label{eq:rk2}
\end{eqnarray}
where $[\cdot]$ represents the jump of a given quantity across the shock. 

Next we expand $u=u_{0}(\xi)+\epsilon u_{1}(\xi,y,t)+O(\epsilon^{2})$,
$v=\epsilon v_{1}(\xi,y,t)+O(\epsilon^{2})$, $s=Dt+\epsilon s_{1}(y,t)+O(\epsilon^{2})$,
where $u_{0}$ is the steady profile, and $D=u_{0s}/2$ is the steady shock speed.
Inserting these expressions into
(\ref{eq:toy-model-2d-saf-1}-\ref{eq:toy-model-2d-saf-2}), and linearizing, we obtain
\begin{alignat}{1}
 u_{1t}+\left(u_{0}-\frac{u_{0s}}{2}\right)u_{1\xi}+u_{0}^{\prime}\?u_{1}+v_{1y} &
    = \left(\frac{\partial f}{\partial u_{s}}\left(x,u_{0s}\right)+
      \frac{1}{2}u_0^{\prime}\right)u_{1}\left(0,t\right),
      \label{eq:2d-linear-problem-1}\\
 u_{1y}-u_{0}^{\prime}\?s_{1y}-v_{1\xi} & =0,\label{eq:2d-linear-problem-2}
\end{alignat}%
 where $u_0^\prime=\frac{d\/u_0}{d\/\xi}$, and the equations apply for
 $\xi \leq 0$ (the uniform state ahead of the shock is unperturbed). The
 linearized shock conditions are
\begin{equation} \label{eq:linearized-rk}
 s_{1t} = \frac{1}{2}u_{1s} \quad\mbox{and}\quad
 s_{1y} = - \frac{v_{1s}}{u_{0s}},
\end{equation}%
 where, as before, the subscript $s$ denotes evaluation at the post-shock state.
It is convenient to introduce
\[
  c_{0} = u_{0} - \frac{1}{2}\?u_{0s},\qquad
  b_{0}=\frac{\partial f}{\partial u_{s}}\left(x,u_{0s}\right)+\frac{1}{2}u_{0}^\prime,
\]
 and then write  (\ref{eq:2d-linear-problem-1}-\ref{eq:2d-linear-problem-2})
 in the form
\begin{align}
 u_{1t}+c_{0}\?u_{1\xi}+u_{0}^{\prime}\?u_{1}+v_{1y} & =b_{0}\?u_{1s},
    \label{eq:2d-linear-problem-2.1}\\
 u_{1y}-v_{1\xi} & =-u_{0\xi}\?v_{1s}/u_{0s}. \label{eq:2d-linear-problem-2.2}
\end{align}
Here $c_{0}$, $b_{0},$ and $u_{0s}$ are determined by the steady-state profile, while
$u_{1s}=u_{1}(0,y,t)$ and $v_{1s}=v_{1}(0,y,t)$ denote the perturbed quantities
evaluated immediately after the shock.

 Linear instability (or stability) is determined by whether or not spatially bounded
 solutions to (\ref{eq:2d-linear-problem-2.1}-\ref{eq:2d-linear-problem-2.2})
 that grow in time exist. In order to answer this question, we use transform
 methods, following the methodology in \cite{Erpenbeck1962}.
 Because (\ref{eq:2d-linear-problem-2.1}-\ref{eq:2d-linear-problem-2.2}) has
 $y$-independent coefficients, we Fourier transform in the transverse direction:
 %
\begin{eqnarray}
 \hat{u}_{t} + c_{0}\?\hat{u}_{\xi} + u_{0}^{\prime}\?\hat{u} + i\?\ell\?\hat{v}
    & = & b_{0}\?\hat{u}_{s}, \label{eq:ft-linear-2d-1}\\
 i\?\ell\?\hat{u} - \hat{v}_{\xi} & = &
    -\frac{u_{0}^{\prime}\?\hat{v}_{s}}{u_{0s}},\label{eq:ft-linear-2d-2}
\end{eqnarray}
where the parameter $\ell$ is the transverse wave number and
\begin{align}
 \hat{u}(\xi,\ell,t) & = \int_{-\infty}^{\infty}e^{-i\ell y}\?u_{1}(\xi,y,t)\?dy,
    \qquad \hat{u}_{s}(\ell,t) = \hat{u}(0,\ell,t),\\
 \hat{v}(\xi,\ell,t) & = \int_{-\infty}^{\infty}e^{-i\ell y}\?v_{1}(\xi,y,t)\?dy,
    \qquad \hat{v}_{s}(\ell,t) = \hat{v}(0,\ell,t).
\end{align}
The Fourier transform of the linearized shock conditions in equation
(\ref{eq:linearized-rk}) is
\begin{equation} \label{eq:ft-linearized-rk}
 \hat{s}_t = \frac{1}{2}\,\hat{u}_s \quad\mbox{and}\quad
 i\?\ell\?u_{0s}\,\hat{s} = - \hat{v}_s
 \quad\mbox{where}\;\;
 \hat{s}(\ell,t) = \int_{-\infty}^{\infty}e^{-i\ell y}\?s_{1}(y,t)\?dy.
\end{equation}
Next, we Laplace transform in $t$ equations
(\ref{eq:ft-linear-2d-1}-\ref{eq:ft-linear-2d-2}), and obtain
\begin{eqnarray}
 \sigma\?U - \hat{u}(\xi,\ell,0) + c_{0}\?U_{\xi} + u_{0}^{\prime}\?U + i\?\ell\?V
    & = & b_{0}\?U_{s},\label{eq:laplace-transformed-ode-1}\\
 i\?\ell\?U - V_{\xi} & = & -\frac{u_{0}^{\prime}\?V_{s}}{u_{0s}},
   \label{eq:laplace-transformed-ode-2}
\end{eqnarray}
where
\begin{align}
 U(\xi,\ell,\sigma) & = \int_{0}^{\infty} e^{-\sigma t}\?\hat{u}(\xi,\ell,t)\?dt, \qquad
    \mbox{with}\;\; U_{s}(\ell,\sigma) = U(0, \ell, \sigma),\\
 V(\xi,\ell,\sigma) & = \int_{0}^{\infty} e^{-\sigma t}\?\hat{v}(\xi,\ell,t)\?dt, \qquad
    \mbox{with}\;\; V_{s}(\ell,\sigma) = V(0, \ell, \sigma).
\end{align}
\begin{note}\label{note:realsigmapos}
 The issue of stability/instability reduces now to the question of the (possible)
 existence of singularities in the Laplace Transform for $\Re(\sigma) > 0$. Thus
 we take $\Re(\sigma) > 0$ in what follows. For that matter: we can calculate
 assuming $\sigma$ real, large, and positive; and then extend the result
 by analytic continuation.\TheoremEnd
\end{note}
Now rewrite (\ref{eq:laplace-transformed-ode-1}-\ref{eq:laplace-transformed-ode-2})
in matrix form
\begin{equation}\label{eq:matrix-form-laplaced-transoformed-ode}
 \mathbf{A} \cdot \frac{d}{d\xi} \mathbf{W} =
    \mathbf{B} \cdot \mathbf{W} + \mathbf{\tilde{F}}(\xi;\sigma,\ell),
\end{equation}
where
\begin{equation}
 \mathbf{W} = \begin{bmatrix}        U\\    V \end{bmatrix},\quad
 \mathbf{A} = \begin{bmatrix} c_{0} & 0\\0 & 1 \end{bmatrix},\quad
 \mathbf{B} = \begin{bmatrix} -\sigma-u_{0}^{\prime} & -i\?\ell\\
                                           i\?\ell & 0 \end{bmatrix},
 \quad\mbox{and}\quad
 \mathbf{\tilde{F}} = \begin{bmatrix}
      \hat{u}(\xi,\ell,0) + b_{0}\?U_{s}\\
      u_{0}^{\prime}\?V_{s}/u_{0s} \end{bmatrix}.
\end{equation}
Two cases arise now, which are different in terms of their analytic complexity.

The first case is the overdriven detonation, where $\zeta>1$ in \eqref{u0s}.
 In this case $c_0 > 0$ for all $\xi<0$. Thus no sonic point exists in the steady
 state, $\mathbf{A}$ is invertible everywhere, and
 (\ref{eq:matrix-form-laplaced-transoformed-ode}) is equivalent to
\begin{equation}\label{eq:overdriven-ode-laplace-matrix}
 \frac{d}{d\xi}\mathbf{W} =
 \mathbf{C} \cdot \mathbf{W} + \mathbf{F}(\xi;\sigma,\ell),
\end{equation}
where
\[
 \mathbf{C} = \mathbf{A}^{-1}\?\mathbf{B} =
 \begin{bmatrix}
   \left(-\sigma-u_{0}^{\prime}\right)/c_{0} & -i\?\ell/c_{0}\\
                                  i\?\ell & 0
 \end{bmatrix} \;\;\mbox{and}\;\;
 \mathbf{F} = \mathbf{A}^{-1}\?\mathbf{\tilde{F}} =
 \begin{bmatrix}
  \left(\hat{u}(\xi,\ell,0) + b_{0}\?U_{s}\right)/c_{0}\\
  u_{0}^{\prime}\?V_{s}/u_{0s}
 \end{bmatrix}.
\]
Here, due to the overdrive assumption, \textbf{$\mathbf{C}$} is a bounded matrix.

The second (harder) case is the Chapman-Jouguet detonation, wherein $c_{0}(\xi)\to 0$
as $\xi \to -\infty$.\footnote{If $f$ is compact support, then the reaction zone has a finite length and $c_0 = 0$ for $\xi$ large enough and negative.  We assume that this is not the case, i.e. the reaction zone is infinite. The results presented can be extended to the cases where $f$ has a compact support.}
\,\,Then, $\mathbf{A}$ is no longer invertible at the sonic point. Hence, a necessary
condition for $\frac{d}{d\xi}\mathbf{W}$ to remain bounded as $\xi\rightarrow-\infty$
is that the right hand side of (\ref{eq:matrix-form-laplaced-transoformed-ode})
become orthogonal to the left eigenvector corresponding to the vanishing
eigenvalue of $\mathbf{A}$, i.e.,
\begin{align*}
 \lim_{\xi\to-\infty} \mathbf{l}_{i} \cdot \mathbf{A}\?\frac{d}{d\xi}\mathbf{W}
   & = \lim_{\xi\to-\infty} \lambda_{i}\?\mathbf{l}_{i} \cdot \frac{d}{d\xi}\mathbf{W}\\
   & = \lim_{\xi\to-\infty} \mathbf{l}_{i} \cdot \left(\mathbf{B} \cdot \mathbf{W} +
       \mathbf{\tilde{F}}(\xi;\sigma,\ell) \right) = 0,
\end{align*}
where $\lambda_{i}$ is the eigenvalue that vanishes as $\xi \to -\infty$, and
$\mathbf{l}_{i}$ is the corresponding eigenvector. Given the form of $\mathbf{A}$,
$\mathbf{B}$, and $\mathbf{F}$, this
 is equivalent to
%
\begin{equation}\label{eq:radiation-cond-2d-toy-model}
 \left(\sigma + u_{0}^{\prime}(\xi)\right)\?U + i\?\ell\?V \rightarrow 0
 \qquad\text{as}\qquad \xi\to-\infty.
\end{equation}%
 This last equation is called the radiation, or boundedness, condition.
 To avoid the difficulties related to the vanishing of $c_{0}$, for
 the remainder of this paper we focus on overdriven detonations.

 Using a fundamental matrix for the homogeneous problem \cite{Coddington1955}, we can
 write the solution to \eqref{overdriven-ode-laplace-matrix} in terms of the boundary
 data at $\xi=0$, and the initial data (encoded into $\mathbf{F}$). That is, 
\begin{equation}\label{eq:LTviaFundSol}
 \mathbf{W}(\xi;\sigma,\ell) = \mathbf{H}(x;\sigma,\ell) \cdot \left[
    \mathbf{H}^{-1}(0;\sigma,\ell) \cdot \mathbf{W}(0;\sigma,\ell) +
    \int_{0}^{x} \mathbf{H}^{-1}(z;\sigma,\ell) \cdot \mathbf{F}(z;\sigma,\ell)
    \?dz\right],
\end{equation}%
 where $\mathbf{H}$ is a non-singular matrix solution of the homogeneous problem:
 $\frac{d}{d\xi}\?\mathbf{H}=\mathbf{C}\cdot\mathbf{H}$ with det$(\mathbf{H})\neq 0$.
 We denote the columns of $\mathbf{H}$ by $\mathbf{h}_1$ and $\mathbf{h}_2$, so that
 $\mathbf{H} = [\mathbf{h}_1, \mathbf{h}_2]$.%

Define $\delta = \lim_{\xi \rightarrow -\infty}\?c_{0}(\xi) > 0$ and
\begin{equation}\label{eq:limxiinf}
 \mathbf{C}_{-\infty} = \lim_{\xi\to-\infty} \mathbf{C}(\xi) =
    \begin{bmatrix}
       -\sigma/\delta & -i\?\ell/\delta\\
              i\?\ell & 0
    \end{bmatrix}.
\end{equation}
$\mathbf{C_\infty}$ is a constant matrix with (generally distinct) eigenvalues
\begin{equation}\label{eq:eigenvalues}
  \lambda_1 = \frac{-\sigma -\sqrt{4\?\delta\?\ell^{2} + \sigma^{2}}}{2\?\delta}
  \quad \mbox{and} \quad
  \lambda_2 = \frac{-\sigma +\sqrt{4\?\delta\?\ell^{2} + \sigma^{2}}}{2\?\delta},
\end{equation}
where we use the principal branch for the square root. Let $\mathbf{q}_1$ and
$\mathbf{q}_2$ be the corresponding eigenvectors. Then, provided that the
limit in (\ref{eq:limxiinf}) is achieved ``fast enough'', $\mathbf{h}_1$ and
$\mathbf{h}_2$ can be selected so that (see Appendix \ref{sec:appendix})
\begin{equation}\label{eq:eigensol_h}
 \mathbf{h}_1 \sim e^{\lambda_{1}\?\xi}\?\mathbf{q}_{1}\quad\mbox{and}\quad
 \mathbf{h}_2 \sim e^{\lambda_{2}\?\xi}\?\mathbf{q}_{2}\quad\mbox{as }\;\xi\to -\infty.
\end{equation}
Since $\Re(\lambda_1) < 0 < \Re(\lambda_2)$, $\mathbf{h}_1$ grows exponentially as
$\xi \to -\infty$, while $\mathbf{h}_2$ decays. Let now $\boldsymbol{\theta}_1$ and
$\boldsymbol{\theta}_2$ be the rows of $\mathbf{H}^{-1}$. Then
\begin{equation}\label{eq:eigensol_theta}
 \boldsymbol{\theta}_j \sim e^{-\lambda_j\?\xi}\boldsymbol{\pi}_j
 \quad\mbox{as }\;\xi\to -\infty,
\end{equation}
where the $\boldsymbol{\pi}_j$ are constant vectors. Furthermore, the $h_j$ and
$\theta_j$ have \emph{analytic dependence} on $\sigma$, for $\Re(\sigma) > 0$  
(Appendix \ref{sec:appendix}).

With the notation above the solution to (\ref{eq:overdriven-ode-laplace-matrix})
can now be written in the form
\begin{equation} \label{eq:LTfinalForm}
 \mathbf{W}(\xi) = \left(\boldsymbol{\theta}_{1}(0) \cdot \mathbf{W}(0) +
    \int_{0}^{\xi} \boldsymbol{\theta}_{1}(z) \cdot \mathbf{F}(s) \? ds\right)
    \mathbf{h}_{1} + \left(\boldsymbol{\theta}_{2}(0) \cdot \mathbf{W}(0) +
    \int_{0}^{\xi} \boldsymbol{\theta}_{2}(z) \cdot \mathbf{F}(s) \? ds\right)
    \mathbf{h}_{2},
\end{equation}
where we only display the dependence on $\xi$. Further, 
boundedness of $\mathbf{W}(\xi)$ as $\xi\rightarrow-\infty$ requires a choice which
eliminates the exponentially growing part of the solution. Specifically:
\begin{equation}\label{eq:bounded-condition-2d-raw}
 \boldsymbol{\theta}_{1}(0) \cdot \mathbf{W}(0) = \lim_{x \rightarrow -\infty}
 \?\int_{x}^{0} \boldsymbol{\theta}_{1}(z) \cdot \mathbf{F}(z) \? dz.
\end{equation}
That this choice is also sufficient to guarantee a bounded solution is less obvious,
but a proof
(given in the context of the reactive Euler equations)
 can be found in \cite{Erpenbeck1962}. Inserting back the definition of
$\mathbf{F}$, we obtain
\begin{equation}\label{eq:bounded-condition-2d}
 \boldsymbol{\theta}_{1}(0) \cdot \mathbf{W}(0) = \int_{-\infty}^{0}
 \boldsymbol{\theta}_{1} \cdot
 \begin{bmatrix} \left(\hat{u}(\xi,\ell,0) + b_{0}\?U_{s}\right)/c_{0}\\
                 u_{0}^{\prime}\?V_{s}/u_{0s}
 \end{bmatrix} dz.
\end{equation}
The Laplace transform of the linearized jump conditions (\ref{eq:ft-linearized-rk})
can be used to relate $V_{s}$ to $U_{s}$ via
\begin{equation}\label{eq:vs-relation-us}
 \sigma \? V_{s} = -i\?\ell\?u_{0s}\left(\frac{1}{2}\?U_{s}+\hat{s}(\ell,0)\right).
\end{equation}%
 Since neither the $h_j$, nor the $\theta_j$ have singular dependence on $\sigma$,
 equation (\ref{eq:LTfinalForm}) shows that singular dependence can enter only
 through the boundary data at $\xi = 0$, that is: $U_s$ and $V_s$. From
 (\ref{eq:vs-relation-us}) it then follows that the stability/instability issue,
 see Note~\ref{note:realsigmapos}, can be decided solely by inspecting the
 (possible) singularities in $U_s = U|_{\xi=0}$.
 To this end, insert now (\ref{eq:vs-relation-us}) into
(\ref{eq:bounded-condition-2d}), and solve for $U_s$. This yields
\begin{equation}\label{eq:FormulaForUs}
 U_s = \frac{\;
    \boldsymbol{\theta}_{1}(0) \cdot
    \begin{bmatrix} 0\\ i\?\ell\?u_{0s}\?\hat{s}(\ell,0) \end{bmatrix}
    + \sigma\int_{-\infty}^{0} \boldsymbol{\theta}_{1} \cdot
    \begin{bmatrix} \hat{u}(\xi,\ell,0)/c_{0} \\
                    i\?\ell\?u_{0}^{\prime}\?\hat{s}(\ell,0)
    \end{bmatrix}\?dz
    \;}{
    \boldsymbol{\theta}_{1}(0) \cdot
    \begin{bmatrix} \sigma \\ -\frac{1}{2}\?i\?\ell\?u_{0s} \end{bmatrix}
    - \int_{-\infty}^{0} \boldsymbol{\theta}_{1} \cdot
    \begin{bmatrix} \sigma\?b_{0}/c_{0} \\ \frac{1}{2}\?i\ell\?u_{0}^{\prime}
    \end{bmatrix}\?dz
    }.
\end{equation}
 Because $\boldsymbol{\theta}_1$  decays exponentially as $\xi\to-\infty$, here both
 the numerator and denominator are regular functions of $\sigma$ for $\Re(\sigma)>0$.
 The possible singularities are thus poles (denominator roots), i.e., solutions to
\begin{equation}\label{eq:bounded-condition-2d-2}
 \boldsymbol{\theta}_{1}(0) \cdot
   \begin{bmatrix} \sigma\\ -\frac{1}{2}\?i\?\ell \end{bmatrix}
   = \int_{-\infty}^{0} \boldsymbol{\theta}_{1} \cdot
   \begin{bmatrix}\sigma\?b_{0}/c_{0}\\ \frac{1}{2}\?i\?\ell\?u_0^{\prime}
   \end{bmatrix}\?dz.
\end{equation}
The function $\boldsymbol{\theta}_{1}$ solves the adjoint homogeneous problem
\begin{eqnarray}
 \frac{d}{d\xi} \boldsymbol{\theta} & = &
    -\mathbf{C}(\xi;\sigma,\ell)^T \cdot \boldsymbol{\theta},
    \label{eq:adjoint-homogenous-problem}\\
 & = & \begin{bmatrix} \left( \sigma + u_{0}^{\prime} \right)/c_{0} & -i\?\ell\\
                                                    i\?\ell/c_{0} & 0
       \end{bmatrix} \boldsymbol{\theta},\nonumber 
\end{eqnarray}
 subject to the condition that $\boldsymbol{\theta}$ should be bounded as
 $\xi \to -\infty$. Thus in this limit $\boldsymbol{\theta}_1$ becomes parallel
 to the eigenvector associated to the eigenvalue with positive real part:
\begin{equation} \label{eq:asymtheta1}
 \boldsymbol{\theta}_1 \sim a\?
 \begin{bmatrix}
    \frac{1}{2}\,\left(\sigma + \sqrt{4\?\ell^2\?\delta + \sigma^2}\right)\\
    i\?\ell
 \end{bmatrix}\?e^{-\lambda_1\?\xi}   \quad\mbox{as}\;\;\xi \to -\infty,
\end{equation}%
 where $a \neq 0$ is a constant. Note that the value of $a$ has no effect on
 equation (\ref{eq:bounded-condition-2d-2}).

When $\ell=0$, the adjoint homogeneous problem (\ref{eq:adjoint-homogenous-problem})
becomes the scalar equation
\begin{equation}
 \frac{d}{d\xi}\mathbf{\theta} = 
    \frac{1}{c_0}\?\left(\sigma + u_0^{\prime}\right)\?\mathbf{\theta},
\end{equation}
with solutions
\begin{equation}\label{eq:homogenous-solution-laplace-1d}
 \mathbf{\theta} = \frac{\theta(0)}{c_0(0)}\?c_0(\xi)\?\exp(-\sigma\?p(\xi)),
 \quad\mbox{where}\;\; p = \int_{\xi}^0 \frac{1}{c_0}\?dz.
\end{equation}
Using (\ref{eq:homogenous-solution-laplace-1d}) in (\ref{eq:bounded-condition-2d-2})
yields
\begin{equation} \label{eq:1D_disp}
  c_{0}(0) =\int_{-\infty}^{0}b_{0}(z)\?e^{-\sigma\?p(z)}\?dz
\end{equation}%
 as the equation for the poles of the Laplace Transform in $\Re(\sigma) > 0$.
 This is the same as the dispersion relation obtained in \cite{Faria2014}, using
 a normal mode approach. It follows that,
in the context of the simple toy model presented here:
\begin{enumerate}
 \item The Laplace transform and normal mode approaches yield the same stability
       criterion in 1D.
 \item The difficulties with the stability analysis for a Chapman-Jouguet detonation, see
  equation (\ref{eq:radiation-cond-2d-toy-model}), do not
  arise in 1D (i.e., $\ell = 0$) for the analog model. The reason is that, in this case,
  the linearized equations can be solved (essentially) explicitly. This analysis,
  using normal modes, can be found in \cite{Faria2014} as well. Surprisingly, the
  answer is the same as for overdriven detonations; i.e., the stability criterion
  is also (\ref{eq:1D_disp}). Thus the question, requiring further investigation is:
  is this still true when $\ell \neq 0$?%
\end{enumerate}
In general, for $\ell \neq 0$, \eqref{adjoint-homogenous-problem} cannot be solved
analytically. Then a numerical method is needed to obtain the eigenvalues (poles of
the Laplace transform), given by the roots $\sigma = \sigma(\ell)$ of
\begin{equation} \label{eq:stabilityFun}
 0 = R(\sigma, \ell) =
 \boldsymbol{\theta}_{1}(0) \cdot
   \begin{bmatrix} \sigma\\ -\frac{1}{2}\?i\?\ell \end{bmatrix}
   - \int_{-\infty}^{0} \boldsymbol{\theta}_{1} \cdot
   \begin{bmatrix}\sigma\?b_{0}/c_{0}\\ \frac{1}{2}\?i\?\ell\?u_0^{\prime}
   \end{bmatrix}\?dz,
\end{equation}
where $\boldsymbol{\theta}_1$ is as in
(\ref{eq:adjoint-homogenous-problem}-\ref{eq:asymtheta1}).
\section{An example}\label{sec:An-example-2d-toy}
 In this section, we select a specific form for the forcing function $f$, and use it 
 to illustrate particular properties of the model. For this purpose, we could use the
 more physically justifiable reaction rate that follows from the approach in
 (\ref{eq:Tunif}-\ref{eq:nonlocal-rection-rate}). However, as pointed out earlier,
 this gives rise to certain complications because we cannot solve
 (\ref{eq:nonlocal-rection-rate}) for $\lambda$ explicitly, which leads to a forcing
 function that is both costly to compute numerically, and hard to understand analytically.
We therefore diverge here from a rational approach and aim instead at simplicity,
choosing $f$ in an \emph{ad hoc }manner.

Motivated by the one dimensional analog model presented in \cite{Kasimov2013,Faria2014},
we choose $f$ to have the form
\[
 f = \frac{q}{\sqrt{4\?\pi\?\beta}\left( 1 + \text{Erf}\left[
        \frac{k}{2\,\sqrt{\beta}}\,\left(\frac{u_s}{u_{0s}}\right)^{-\alpha}
        \right]\right)}\?\exp\left[ - \frac{1}{4\,\beta}\?\left( x-x_s+ k\?\left(
        \frac{u_{0s}}{u_s} \right)^{\alpha} \right)^2\right]
\]
 for $x < x_s$, where $q > 0$, $k > 0$, $\beta > 0$, and $\alpha > 0$ are
 parameters, and $f = 0$ ahead of the shock, $x > x_s$.

 The form above mimics important features of the reaction rate, controlled by
 the parameters. It has an induction length (determined by $k$), a reaction-zone
 width (governed by $\beta$), and a heat release $q$ (note that the total energy release is constant, 
 $\int f\/dx = q/2$,  as in \cite{Faria2014}). Finally, $\alpha$ regulates the sensitivity to
 the shock velocity (and plays a role similar to that of the activation
 energy).

It is convenient to rescale the variables as follows 
\[
 u = u_{0s}\?\tilde{u},\quad v = u_{0s}^{3/2}\?\tilde{v},\quad x = k\?\tilde{x},
 \quad y = k\?\tilde{y}/\sqrt{u_{0s}},\quad t = k\?\tilde{t}/u_{0s},
\]
where $u_{0s}$ is given by \eqref{u0s}. Then the model equations take the form
\begin{eqnarray}
 \tilde{u}_{\tilde{t}} + \frac{1}{2} \left(\tilde{u}^{2}\right)_{\tilde{x}} +
   \tilde{v}_{\tilde{y}} & = & \tilde{f}(\tilde{x}-\tilde{x}_{s}, \tilde{u}_{s}),
   \label{eq:toy-model-2d-1-1}\\
 \tilde{v}_{\tilde{x}} - \tilde{u}_{\tilde{y}} & = & 0, \label{eq:toy-model-2d-2-1}
\end{eqnarray}
where the rescaled forcing function $\tilde{f}$ is
\begin{equation}\label{eq:forcing_dimensionless-1}
 \tilde{f}(\tilde{x} - \tilde{x}_{s}, \tilde{u_{s}}) =
   \frac{1}{4\zeta^{2}\left(1+\text{Erf}\left[
      \frac{\tilde{u}(0,\tilde{t})^{-\alpha}}{2\sqrt{\tilde{\beta}}}\right]\right)}
   \frac{1}{\sqrt{4\?\pi\?\tilde{\beta}}}
   \?\exp\left[-\frac{1}{4\?\tilde{\beta}}\?\left(
      \tilde{x} - \tilde{x}_{s} +
      \left( \tilde{u}\left(0,\tilde{t}\right) \right)^{-\alpha}\right)^{2}\right],
\end{equation}
$\tilde{\beta} = \beta/k$, and $\zeta$ is the overdrive parameter introduced
in (\ref{eq:u0s}).

 Therefore, the wave dynamics is controlled by three parameters:  $\alpha$
 -- sensitivity of the reaction rate to variations in $u$ at the shock,
 $\tilde{\beta}$ -- ratio of the reaction zone length to the induction zone length,
 and $\zeta$ -- the degree of overdrive.
From here on we use the dimensionless variables,
and drop the tilde notation.

\begin{figure}[htb!]\vspace*{-0.8em}
 \centering{}
 \includegraphics[width=4in,height=2in]{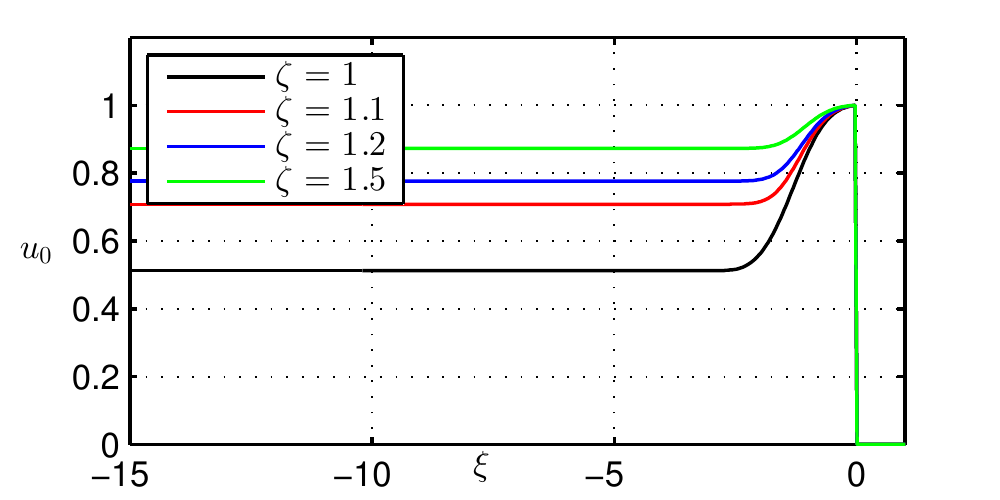}
 \protect\caption{\label{fig:Steady-state-solution-profiles-overdrive}Steady-state
 solution profiles for (\ref{eq:toy-model-2d-1}-\ref{eq:toy-model-2d-2})
 as the overdrive is varied, while keeping all other parameters fixed. }
\end{figure}

 In the dimensionless equations, the overdrive $\zeta$ scales the amplitude of the
 source term. Thus, the influence of chemical reactions vanishes as
 $\zeta \to \infty$, and the wave approaches an inert shock. This is illustrated
 by \figref{Steady-state-solution-profiles-overdrive}, showing the effect of the
 overdrive factor on the steady detonation profile. For large enough overdrive,
 $u$ becomes nearly constant behind the precursor shock, and the wave is primarily
 sustained by the imposed left boundary condition.

In the next subsection we focus on the role played by multi-dimensional effects on:
 (1) the linear stability properties, and
 (2) the full nonlinear dynamics of the solutions to
     (\ref{eq:toy-model-2d-1-1}-\ref{eq:toy-model-2d-2-1}).
The goal is to show that multi-dimensional instabilities typically dominate
one-dimensional instabilities, and that these instabilities lead to the formation
of complex patterns in the solutions to
 (\ref{eq:toy-model-2d-1-1}-\ref{eq:toy-model-2d-2-1}).
We start with the linear stability analysis.
\subsection{Multidimensional linear stability analysis}\label{subsec:MultiDstabAn}
\begin{figure}[htb!]\vspace*{-0.8em}
\noindent\centering{}
\includegraphics[width=2.8in,height=2in]{%
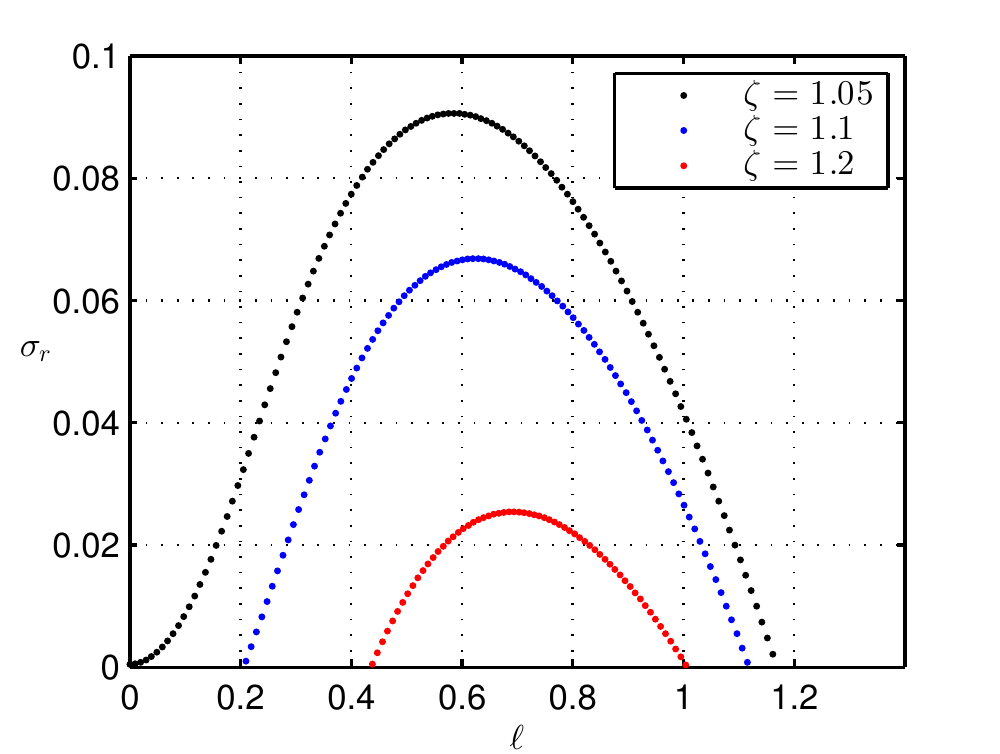}\;\;
\includegraphics[width=2.8in,height=2in]{%
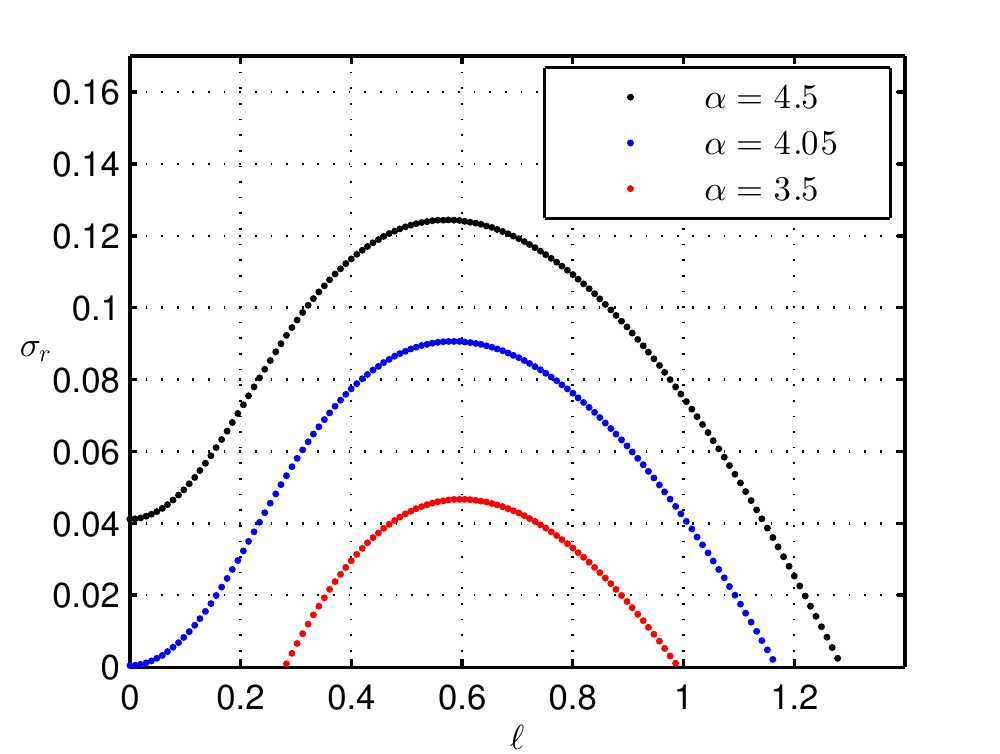}
\protect\caption{\label{fig:Dispersion-relation-for-vary-overdrive-alpha}
Growth rate, $\sigma_r$, as a function of the transverse wavenumber, $l$\,: (a) at $\beta=0.1$, $\alpha=4.05$, and varying degree of overdrive, $\zeta$; 
(b) at $\beta=0.1$, $\zeta=1.05$, and varying $\alpha$.}
\end{figure}

As shown in Section \ref{sub:Linear-stability-2d-toy},
 the stability question for the steady state solutions of
 (\ref{eq:toy-model-2d-1-1}-\ref{eq:toy-model-2d-2-1}) is decided by the roots of
 the stability function $R(\sigma\/,\,\ell)$, defined in (\ref{eq:stabilityFun}).
 If zeros exist with $\Re(\sigma) > 0$ for any admissible $\ell$ (discrete or continuous depending on whether the domain in $y$ is bounded or not), then the profile is unstable.
 Otherwise, it is stable.
 The main difficulty with finding the roots to (\ref{eq:stabilityFun}) stems from
 the fact that $\theta_1$, defined in
 (\ref{eq:adjoint-homogenous-problem}-\ref{eq:asymtheta1}), must be computed
 numerically. Thus, each evaluation of $R(\sigma,\ell)$ requires solving the linear
 ODE system (\ref{eq:adjoint-homogenous-problem}) to obtain $\theta_1$, followed by
 a numerical evaluation of the integral in the definition of $R$. This makes
 parametric studies for varying $\alpha$ and $\beta$ computationally costly.
 We start by investigating the overdrive effect on the wave stability. As pointed
 out earlier, because detonations become ``closer'' to inert shocks the larger the
 overdrive is, we expect the overdrive to have a stabilizing effect. This is
 confirmed by \figref{Dispersion-relation-for-vary-overdrive-alpha}(a), where the
 growth rate $\sigma_{r} = \Re(\sigma(\ell))$, as a function of the transverse wave number $\ell$,
 for $\beta=0.1$ and $\alpha=4.05$, is plotted for increasing values of the overdrive
 parameter $\zeta=1.05$, $1.10$, and $1.20$. For the calculations in this subsection, it turns out that
    equation (\ref{eq:stabilityFun}) has a single root in $\Re(\sigma) > 0$. In
    general, the growth rate is the maximum value over all the roots of
    $\Re(\sigma)$.
 It is seen in \figref{Dispersion-relation-for-vary-overdrive-alpha}(a) that the maximum growth rate occurs at a value $\ell \neq 0$.
 In particular, waves that are stable to purely longitudinal disturbances
 $\left(\ell=0\right)$, can be unstable to 2D perturbations. Thus, generally, we expect
 multi-dimensional effects to dominate over 1D effects.

We also study the effect of $\alpha$ on the stability of the waves. Since $\alpha$
measures the sensitivity of the forcing to changes in the lead shock strength,
we expect that larger values of $\alpha$ will augment the instability. Specifically, 
the growth rate $\sigma_r$ should increase with $\alpha$. This is precisely what is
observed in \figref{Dispersion-relation-for-vary-overdrive-alpha}(b).
This figure also shows that $\alpha$ seems to have very little effect on the value
of $\ell$ corresponding to the most unstable transverse mode. In particular, for
the parameters in \figref{Dispersion-relation-for-vary-overdrive-alpha}(b), the
most unstable mode occurs for $\ell\approx0.6$, regardless of the value $\alpha$
takes. This behavior is similar to the observation in \cite{Faria2014} for the
1D stability. That is, $\alpha$ has very little effect on the imaginary part
of the unstable eigenvalues (recall that \cite{Faria2014} uses a mode based
stability analysis).
\subsection{Numerical simulations}\label{subsec:NumSim}
The linear stability results presented in \S~\ref {subsec:MultiDstabAn} suggest that 2D
effects play an important role in the solutions of
(\ref{eq:toy-model-2d-1-1}-\ref{eq:toy-model-2d-2-1}). Here, we investigate what
happens after the onset of the instabilities. We use numerical simulations
to investigate the large time limit of solutions to
(\ref{eq:toy-model-2d-1}-\ref{eq:toy-model-2d-2}), which start from initial
conditions given by (linearly unstable) ZND profiles
(\ref{sub:Traveling-wave-solutions-2d-toy}). These calculations show that the
two-dimensional analog model exhibits many of the interesting features of real
multi-dimensional detonations.

\begin{figure}[htb!]\vspace*{-0.8em}
\noindent \centering{}
\subfloat[$u$ at time $t=800$.]{\includegraphics[width=2.8in,height=1.8in]{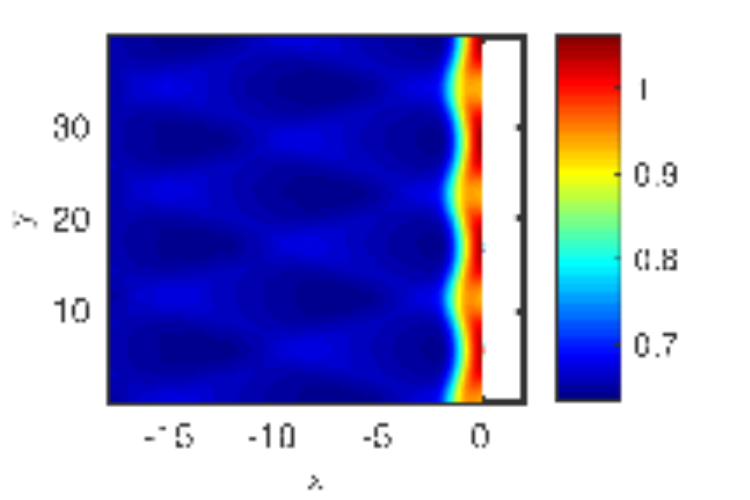}} 
\;\;
\subfloat[$f$ at time $t=800$.]{\includegraphics[width=2.8in,height=1.8in]{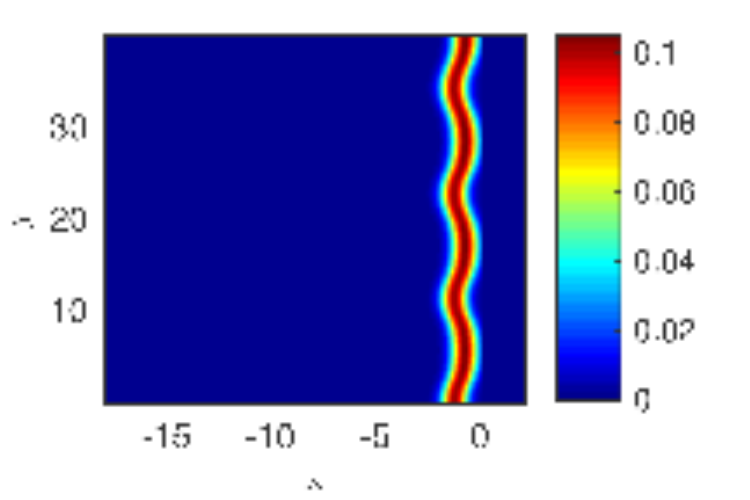}} 
\protect\caption{\label{fig:regular-cells-2d-toy}%
Large time dynamics for $\beta=0.1$, $\alpha=3.5$, and $\zeta=1.05$. This
corresponds to a weakly unstable regime, with regular structures forming,
and small transverse velocities.}
\end{figure}

The numerical simulation of (\ref{eq:toy-model-2d-1}-\ref{eq:toy-model-2d-2}) poses
some problems. On the one hand, because it is a genuinely nonlinear hyperbolic system,
shocks can form even when starting from smooth initial data. On the other hand, it has
characteristic surfaces which are orthogonal to time. Thus the evolution in time is a
nonlocal process. A scheme capable of dealing with shocks is needed, but it
cannot be fully explicit since it is impossible to satisfy a typical CFL condition.
Even with no reaction, $f\equiv0$, when the equations reduce to
\begin{eqnarray}
 u_t+\left(\frac{1}{2}\,u^{2}\right)_{x}+v_{y} & = & 0,\label{eq:dKP1}\\
 v_{x}-u_{y} & = & 0,\label{eq:dKP2}
\end{eqnarray}
developing an appropriate algorithm is not straightforward (e.g., see
\cite{Hunter2000,tabak1994focusing}). We solve
(\ref{eq:toy-model-2d-1}-\ref{eq:toy-model-2d-2}) using an algorithm
based on a semi-implicit time discretization.
The approach is discussed in \cite{Faria2014b}, where a detailed numerical study of
asymptotic equations similar to (\ref{eq:toy-model-2d-1}-\ref{eq:toy-model-2d-2})
is performed. Finally, we note that using a shock-fitting algorithm is not as simple
as in the one-dimensional case \cite{Faria2014}, due to the presence of transverse
shocks. We therefore adopt a shock capturing approach.
%
Equations (\ref{eq:toy-model-2d-1}-\ref{eq:toy-model-2d-2}) are solved in an
inertial frame of reference moving with constant speed $D=1/2$, which is the speed
of the unperturbed ZND wave. The top and bottom ($y$ coordinate) boundary conditions
are that of a rigid wall. An inflow boundary condition on the right, and an outflow
on the left are given.
We have verified that, when linear stability predicts a stable traveling wave (e.g.,
when $\alpha$ is small enough), the numerical solver is able to correctly capture
the ZND structure and wave speed, even when a small amplitude perturbation is
initially imposed on the exact ZND profile.
Interesting dynamics are observed for parameters selected so that the initial
data is unstable.
If the ZND wave is only weakly unstable (meaning the parameters are close to the
neutral stability boundary), very regular multi-dimensional patterns are observed:
see \figref{regular-cells-2d-toy}. 
Figure~\ref{fig:regular-cells-2d-toy}(b)
displays regions  where the induction zone length (distance between the lead
shock and the peak of $f$) is significantly reduced relative to the ZND theory.
In these regions the energy is released shortly after the lead shock. The
associated transverse waves, however, appear to be smooth:
\figref{regular-cells-2d-toy}(a).

\begin{figure}[htb!]\vspace*{-0.8em}
\subfloat[$u$ at $t=1000$, for $\zeta=1.10$.]{%
\includegraphics[width=2.8in,height=1.8in]{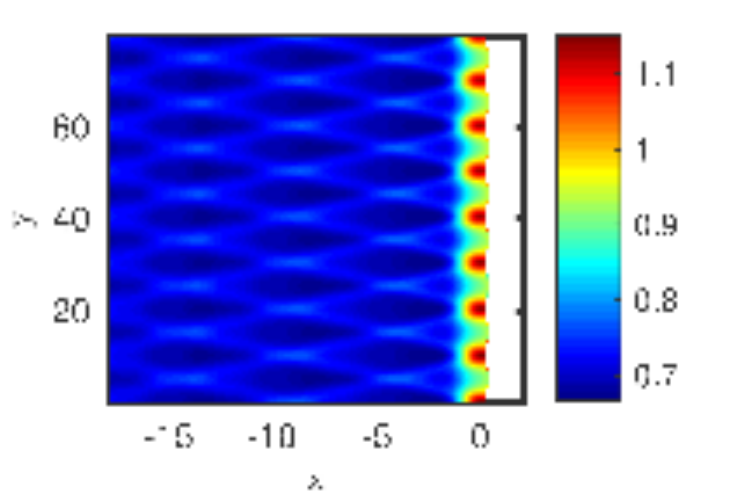}} 
\subfloat[$f$ at $t=1000$, for $\zeta=1.10$.]{%
\includegraphics[width=2.8in,height=1.8in]{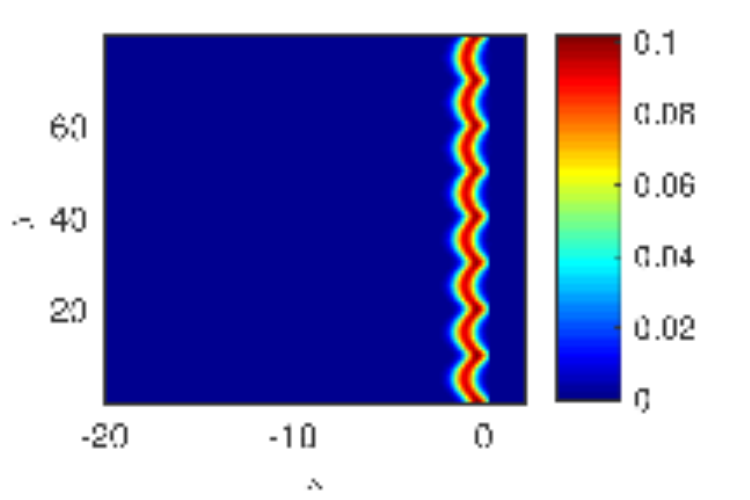}}\\ 
\subfloat[$u$ at $t=1000$, for $\zeta=1.05$.]{%
\includegraphics[width=2.8in,height=1.8in]{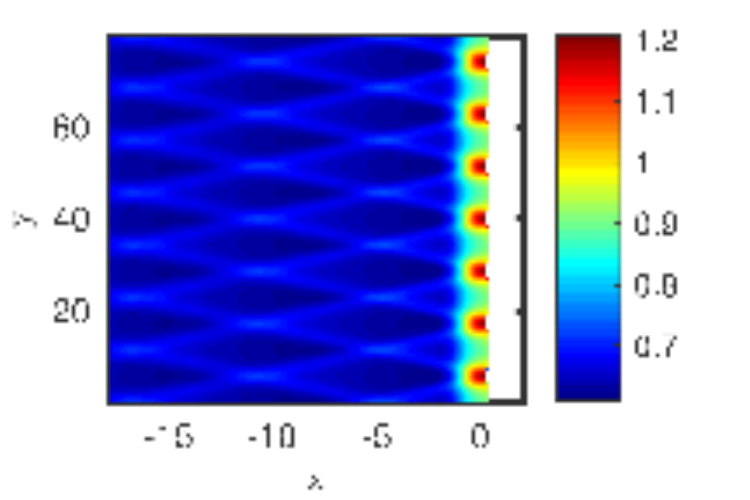}} 
\subfloat[$f$ at $t=1000$, for $\zeta=1.05$.]{%
\includegraphics[width=2.8in,height=1.8in]{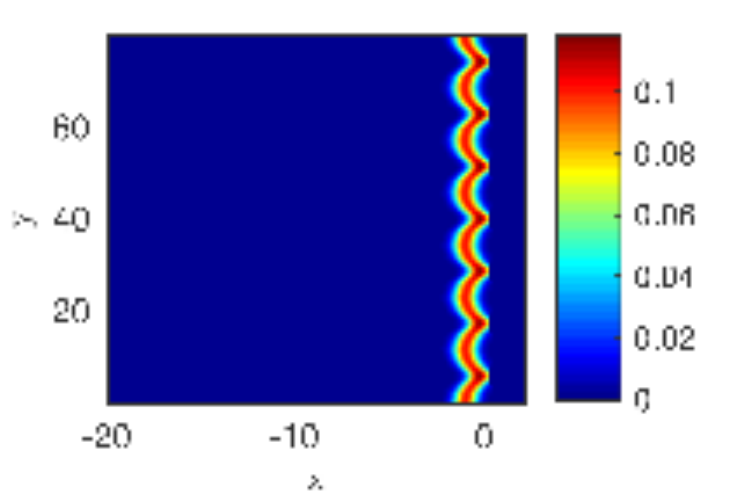}} 
\protect\caption{\label{fig:detonation-vary-overdrive}%
Cellular patterns at varying overdrive, $\zeta$, and $\alpha=4.05$, $\beta=0.1$.
}
\end{figure}

Next, we investigate the effect of the overdrive parameter $\zeta$. Waves closer to
the Chapman-Jouguet case are more unstable and display stronger transverse
variations. Furthermore, for smaller overdrives, the cells in the patterns are larger.
These findings are consistent with the linear stability prediction in
\figref{Dispersion-relation-for-vary-overdrive-alpha}. There it can be seen that:
(1) smaller $\zeta$ corresponds to larger growth rates, and (2) the most unstable
wavelength increases with decreasing degree of overdrive. 
Figure~\ref{fig:detonation-vary-overdrive} illustrates these points.

\begin{figure}[htb!]\vspace*{-0.8em}
\subfloat[$u$ at $t=1000$, for $\alpha=4.05$.]{%
\includegraphics[width=2.8in,height=1.8in]{%
Figures/alpha_4\lyxdot 05_u_t_1000_y_80}} 
\subfloat[$f$ at $t=1000$, for $\alpha=4.05$.]{%
\includegraphics[width=2.8in,height=1.8in]{%
Figures/alpha_4\lyxdot 05_omega_t_1000_y_80}}\\
\subfloat[$u$ at $t=1000$, for $\alpha=4.30$.]{%
\includegraphics[width=2.8in,height=1.8in]{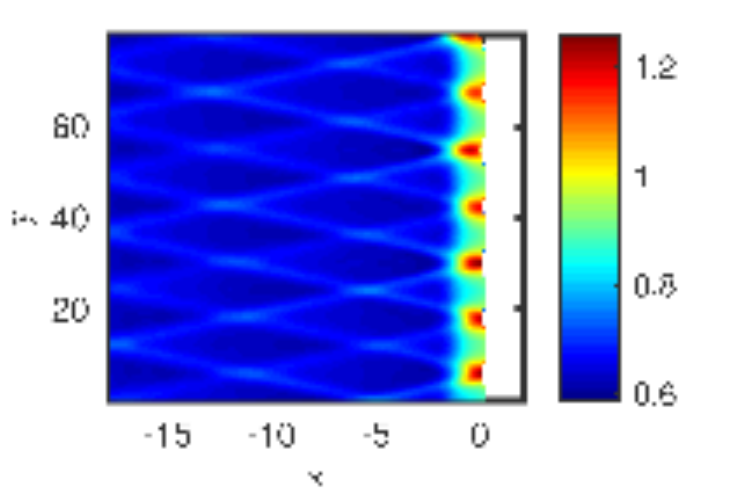}}
\subfloat[$f$ at $t=1000$, for $\alpha=4.30$.]{%
\includegraphics[width=2.8in,height=1.8in]{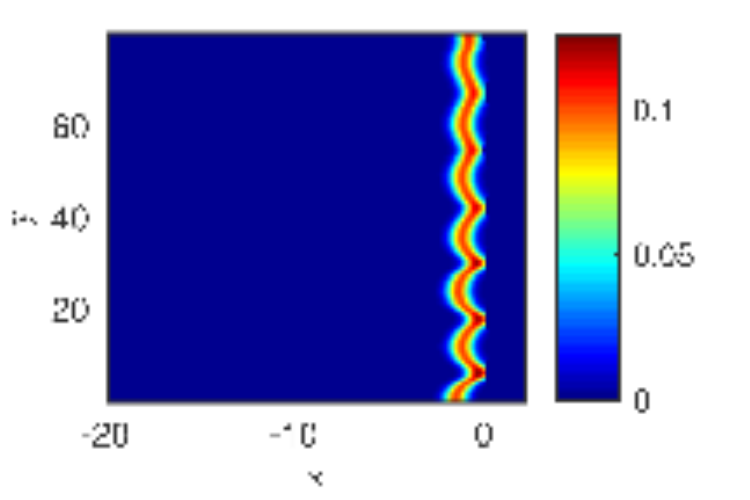}}\\
\subfloat[$u$ at $t=1000$, for $\alpha=4.50$.]{%
\includegraphics[width=2.8in,height=1.8in]{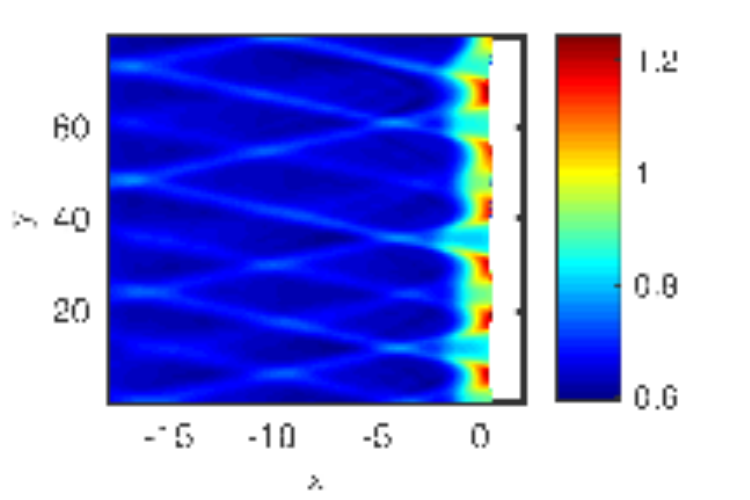}}
\subfloat[$f$ at $t=1000$, for $\alpha=4.50$.]{%
\includegraphics[width=2.8in,height=1.8in]{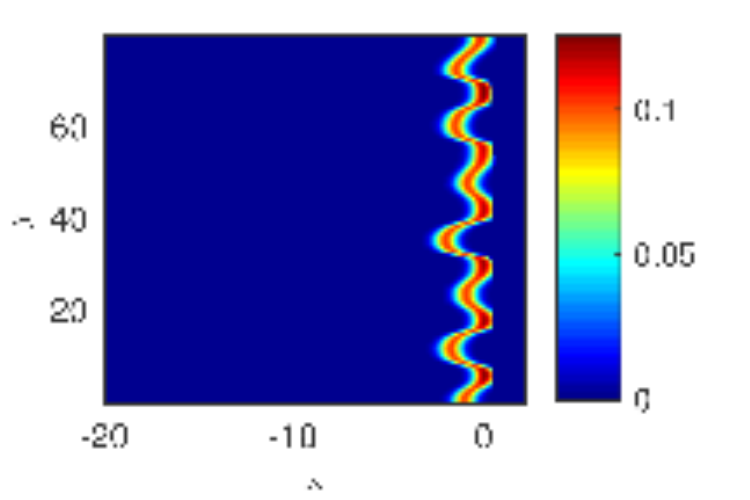}}
\protect\caption{\label{fig:detonation-overdrive-vary-alpha}%
Large time behavior of the solutions for $\zeta=1.05$, $\beta=0.1$,
and varying $\alpha$.}
\end{figure}

 Figure~\ref{fig:detonation-overdrive-vary-alpha} explores the effect of $\alpha$,
 the shock-state sensitivity of the reaction rate, on the detonation wave structure
 and stability. As in the one dimensional case \cite{Faria2014}, larger values of
 $\alpha$ are associated with more complex dynamics. As $\alpha$ grows, so does the
 strength of the patterns generated, which also become more irregular.
 Figure~\ref{fig:detonation-overdrive-vary-alpha}(e-f) shows an example of an
 irregular cellular detonation wave.

These figures, qualitatively, match well the cellular patterns observed for gaseous detonations in dilute mixtures \cite{Lee2008,Fickett2012,OranBoris}. However, a quantitative characterization of the two-dimensional dynamics, of the sort found for the 1D case in \cite{Faria2014}, is challenging. As $\alpha$ increases, the solutions are seen to go through a series of bifurcations, and transition: from (i) a very regular, periodic, structure; to (iii) irregular and
 aperiodic behavior; through (ii) intermediate stages of increasing complexity.
 These stages are illustrated by each of the rows in
 \figref{detonation-overdrive-vary-alpha}. Panels (a-b) correspond to a very simple,
 periodic, structure. Panels (c-d) show some extra structure, easily detectable in
 the $y$-direction in panel (d). Finally, the waves in panels (e-f) are fairly
 irregular.

\begin{figure}[htb!]\vspace*{-1.0em} 
\begin{centering}
\subfloat[]{\includegraphics[width=5.0in,height=1.5in]{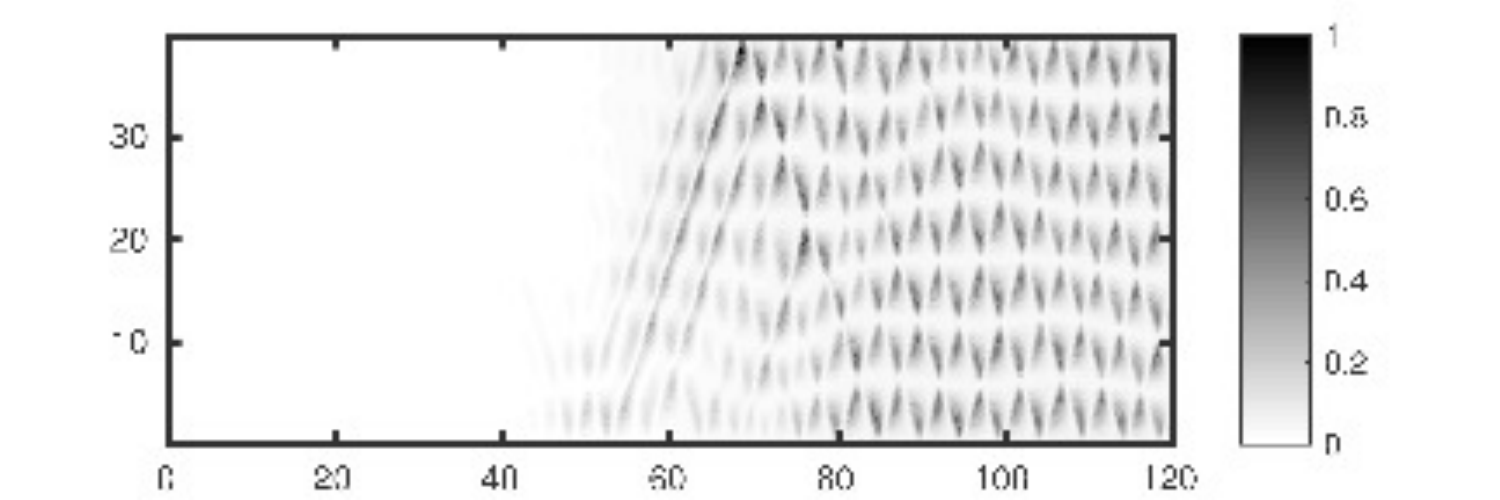}} 
\par\end{centering}
\begin{centering}
\subfloat[]{\includegraphics[width=5.0in,height=1.5in]{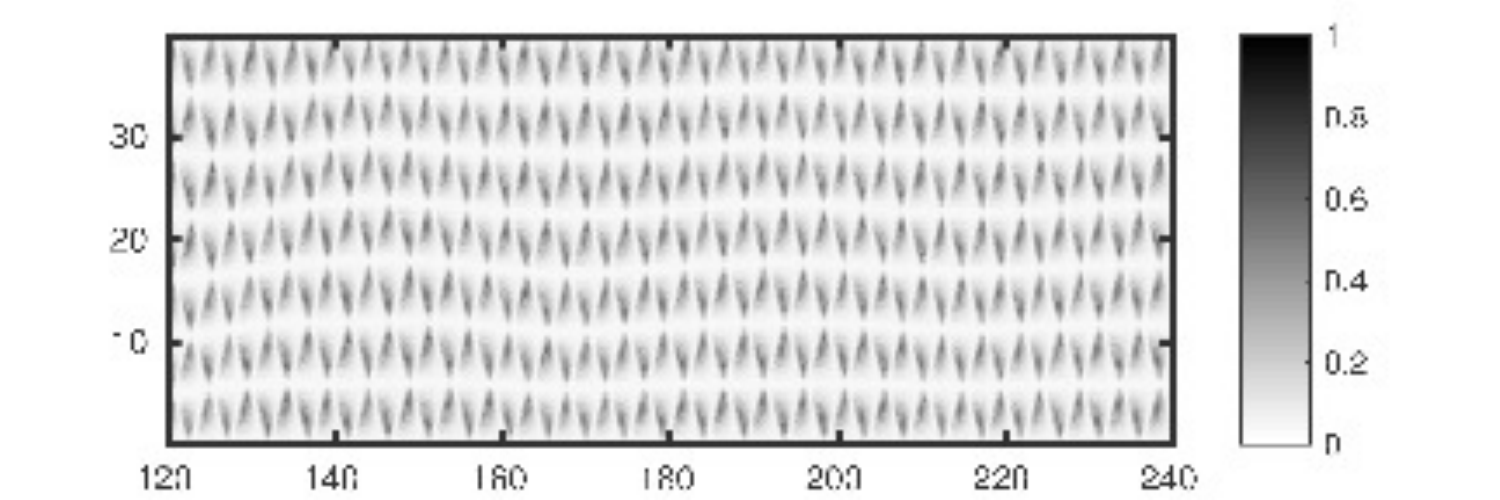}} 
\par\end{centering}
\begin{centering}
\subfloat[]{\includegraphics[width=5.0in,height=1.5in]{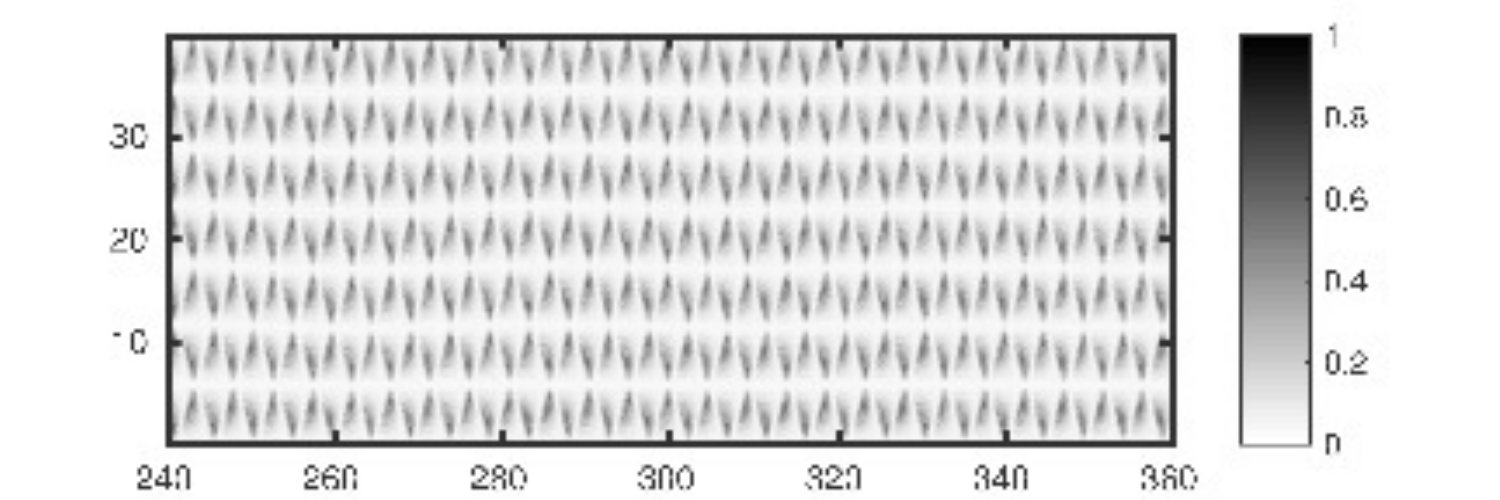}} 
\par\end{centering}
\centering{}\protect\caption{\label{fig:v-trace-1}%
Trace of $|\nabla v|$, as in \eqref{tracegradv}, for $\alpha=4.1$, with
$\beta=0.1$ and $\zeta=1.05$. The horizontal and vertical plot axes are
$x_l$ and $y_l$, and the intensity gray bar is normalized to $[0\/,\,1]$.
Each panel corresponds to a third of the plot range $0 \leq x_l \leq 360$.}
\end{figure}

\begin{figure}[htb!]\vspace*{-0.8em} 
\begin{centering}
\subfloat[]{\includegraphics[width=5.0in,height=1.5in]{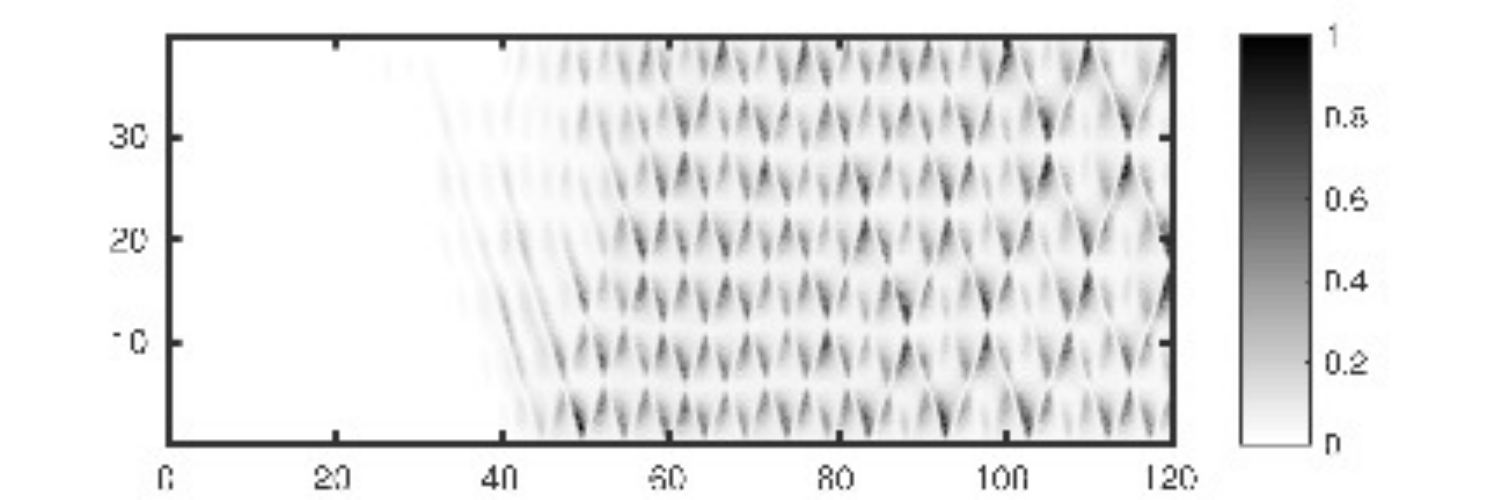}} 
\par\end{centering}
\begin{centering}
\subfloat[]{\includegraphics[width=5.0in,height=1.5in]{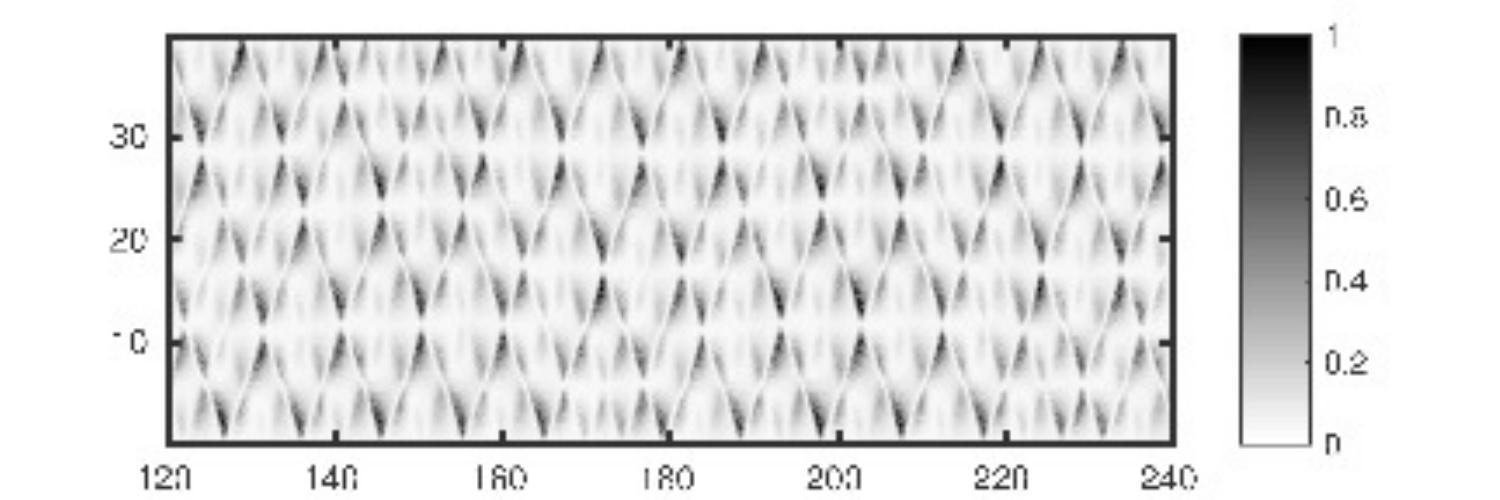}} 
\par\end{centering}
\begin{centering}
\subfloat[]{\includegraphics[width=5.0in,height=1.5in]{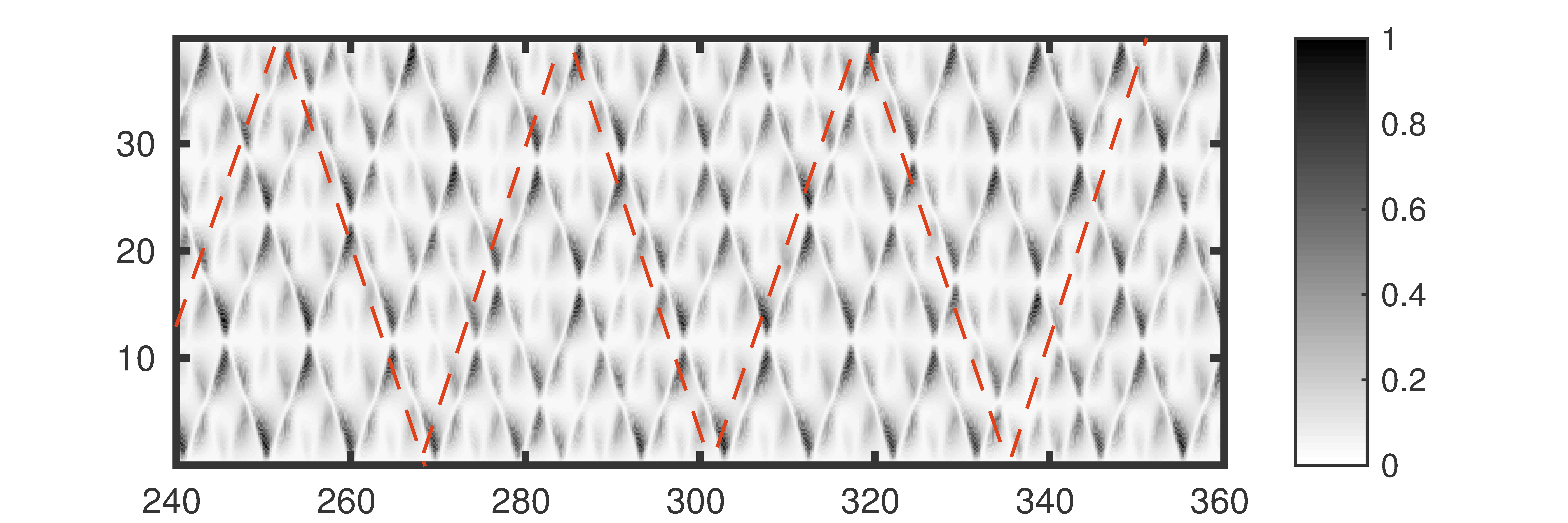}} 
\par\end{centering}
\centering{}\protect\caption{\label{fig:v-trace-2}%
Trace of $|\nabla v|$, as in \eqref{tracegradv}, for $\alpha=4.5$, with
$\beta=0.1$ and $\zeta=1.05$. The horizontal and vertical plot axes are
$x_l$ and $y_l$, and the intensity gray bar is normalized to $[0\/,\,1]$.
Each panel corresponds to a third of the plot range $0 \leq x_l \leq 360$.
The dashed line in panel (c) tracks a transverse wave as it bounces
back-and-forth from the boundaries.
}
\end{figure}

\begin{figure}[htb!]\vspace*{-0.8em} 
\begin{centering}
\subfloat[]{\includegraphics[width=5.0in,height=1.5in]{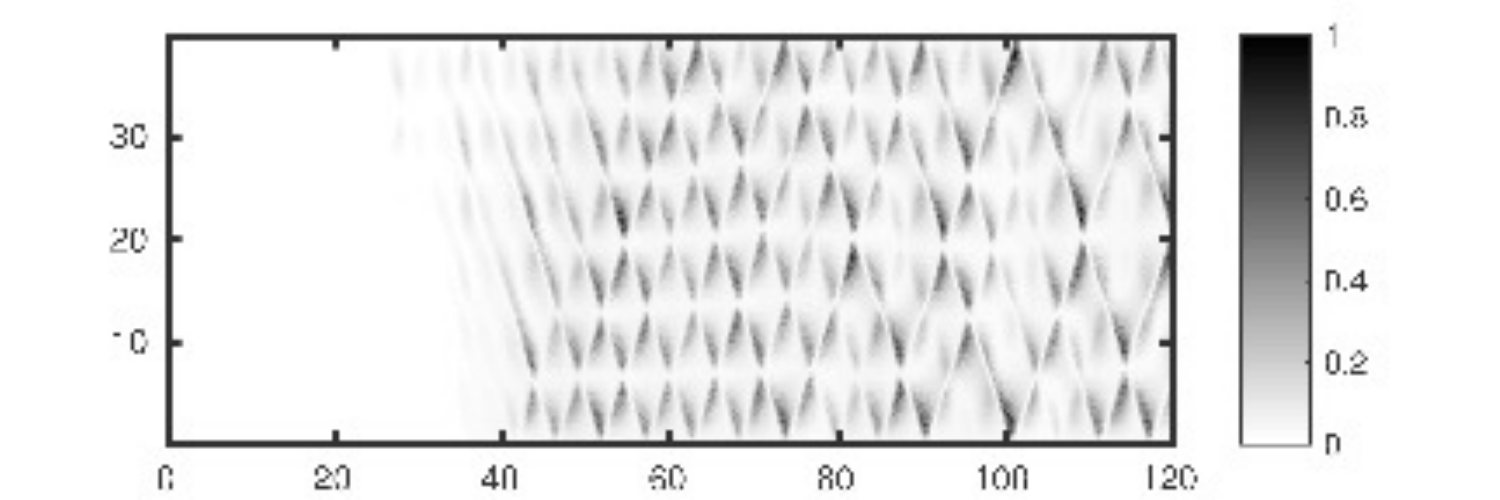}} 
\par\end{centering}
\begin{centering}
\subfloat[]{\includegraphics[width=5.0in,height=1.5in]{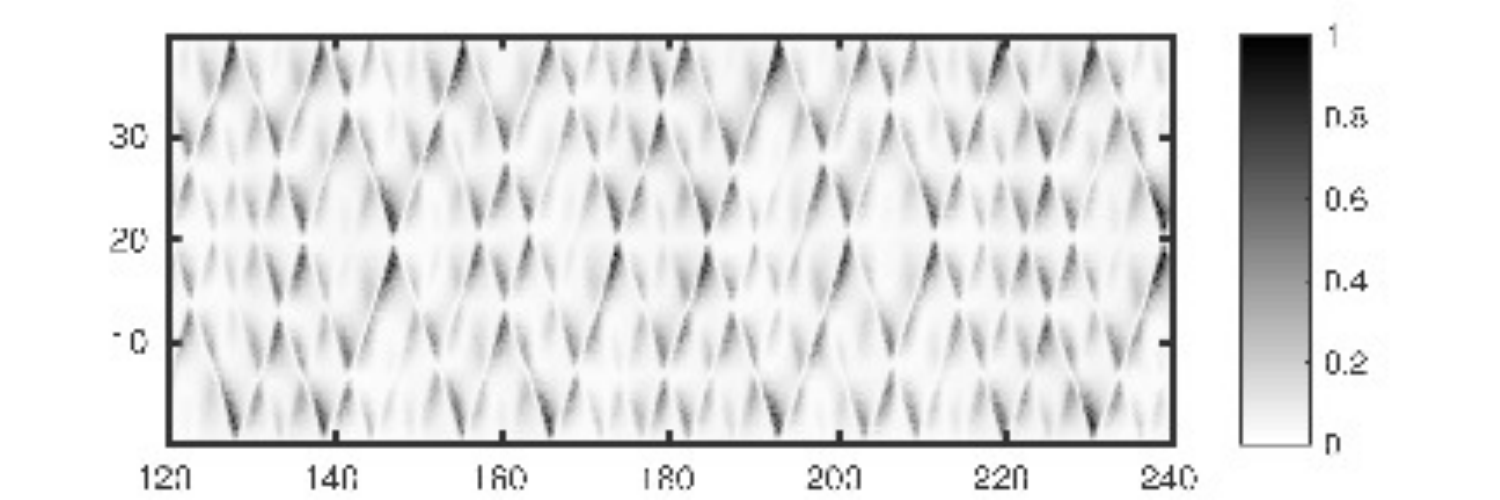}} 
\par\end{centering}
\begin{centering}
\subfloat[]{\includegraphics[width=5.0in,height=1.5in]{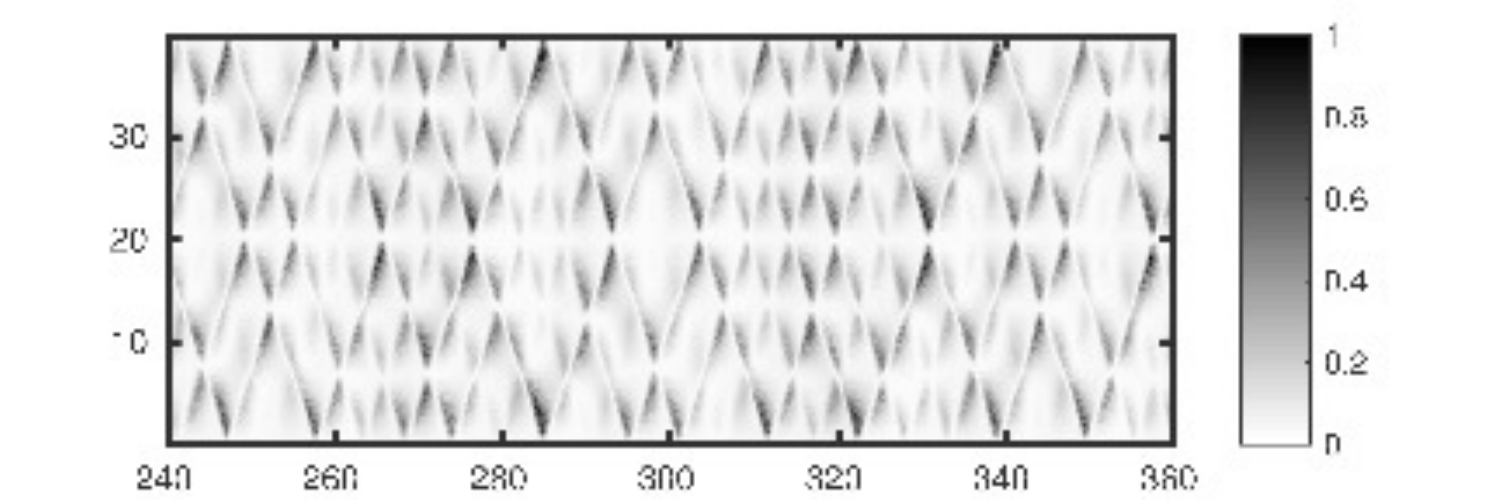}} 
\par\end{centering}
\centering{}\protect\caption{\label{fig:v-trace-3}%
Trace of $|\nabla v|$, as in \eqref{tracegradv}, for $\alpha=4.8$, with
$\beta=0.1$ and $\zeta=1.05$. The horizontal and vertical plot axes are
$x_l$ and $y_l$, and the intensity gray bar is normalized to $[0\/,\,1]$.
Each panel corresponds to a third of the plot range $0 \leq x_l \leq 360$.}
\end{figure}

 Figures~\ref{fig:v-trace-1} through \ref{fig:v-trace-3} display the time
 history of the gradient of $v$, for successively higher values of $\alpha$. These plots are produced as follows.
 First, we convert the solutions from the ZND-wave moving frame (used for the
 numerical calculation) back to the ``lab'' frame $(x_l\/,\,y_l)$. The calculation is
 run for a long enough time, $0\leq t\leq T_f$, to insure that transient effects are
 gone. This results into the lead shock moving (in the lab frame) several hundred
 units: the shock starts at $x_l \approx 20$, and it is beyond $x_l = 400$ by $t=T_f$. 
 Second, we construct the ``trace'' of $|\nabla v|$, using the formula
 \begin{equation}\label{eq:tracegradv}
   Tr(x_l\/,\,y_l) = \int_{0}^{T_f} |\nabla v(x_l\/,\,y_l\/,\,t)|\?dt.
 \end{equation}
This formula makes sense even if $v$ is discontinuous, and is equivalent to the total variation in time at a given point $(x_l,y_l)$. The trace is then plotted for $x_l$ in the interval $0 \leq x_l \leq 360$.
 The results are displayed in figures~\ref{fig:v-trace-1}-\ref{fig:v-trace-3}.
 Note that, since the transverse waves are stronger closer to the lead shock (and
 decay away from it), moving left to right in these plots is (roughly) equivalent
 to exploring the evolution in time of the transverse waves. 
 The process described above is a numerical analog of the experimental soot
 foil records of detonations. A related approach was used in \cite{Oran1998}, where the authors traced the time integral of a specific energy release at a given location.
 
 The plots illustrate that the regularity of the final cellular structure produced
 depends on the distance to the stability boundary (larger the larger $\alpha$ is),
 as well as several other properties.
 For $\alpha=4.1$, \figref{v-trace-1} shows that the cells become regular after a
 sufficiently long time. Initially (top panel), cells of differing strengths coexist,
 and there is some irregularity in the pattern. After a while, the pattern becomes
 nearly regular, with only a low frequency modulation remaining. This is apparent
 in the middle panel, which exhibits slight horizontal waviness (easier to see along
 the whitest regions of the plot). This transient decays slowly, and eventually
 disappears --- the bottom panel shows very little of it.

 Figure~\ref{fig:v-trace-2}, corresponding to $\alpha=4.5$, exhibits an interesting
 phenomenon. Initially (top panel), the cells are relatively small. However, after a
 while, cells that are about twice as large as the original ones appear (middle
 panel), while the smaller ones gradually decay. Eventually (bottom panel), the
 larger cells completely take over the pattern. This behavior is similar to the one
 found in detonation simulations with the reactive Euler equations. That is: the
 initial cell sizes are found to be consistent with the length scale provided by
 the most unstable linear mode, but the long time behavior is dominated by cells
 that can be significantly larger \cite{gamezo1999formation,sharpe2000two}. 
  In fact, to further stress this point, note that: in each of the figures the
 initial cell size is, roughly, the same. This corresponds to the fact that
 $\ell$ for the most unstable linear mode does not change much when $\alpha$
 changes --- see figure~\ref{fig:Dispersion-relation-for-vary-overdrive-alpha}.
 On the other hand, the onset of the cells happens earlier as $\alpha$ grows,
 because the growth rate of the most unstable linear mode also increases with $\alpha$.

When $\alpha$ is increased further, as in \figref{v-trace-3}, cell regularity is
lost (at least for the time scales our computations explored).
Here, again, the average cell size increases after the initial transient stage.
%

 Another peculiar feature, which can be seen in \figref{v-trace-2}, is that the
 amplitude of the transverse waves has large variations as they bounce
 back-and-forth from the channel walls. This can be verified by tracing the path
 of any single wave as it bounces from the walls (see the dashed line
 in the third panel of figure~\ref{fig:v-trace-2}). It is possible that this is
 simply another transient phenomenon, decaying after a very long time. However,
 this may be a part of the multi-modal dynamics of the wave, being
 produced by the collisions between transverse waves and/or the interaction of the longitudinal and transverse instabilities, which are both present in this particular case.
%
\section{Conclusions}
A simple two-dimensional analog model for detonations was introduced. The model
consists of a non-locally forced Burgers-like equation, coupled with a
zero-vorticity equation, and extends our earlier model for one-dimensional pulsating
detonations. The model was shown to be capable of capturing the multi-dimensional
principal characteristics of cellular detonations, wherein transverse shocks and
associated triple-point formation play an important role.

A Laplace transform approach was employed to perform a linear stability analysis of
the model traveling wave solutions. It was shown that the dispersion relation for
the unstable modes can be written as an integral equation, much like the explicit
1D result derived in \cite{Faria2014}. However, in 2D evaluation of the dispersion
relation requires a computationally intensive procedure, as it is also the case for
detonations modeled by the full reactive Euler equations.
In particular: the overdrive factor was shown to have a stabilizing effect (same as
for the reactive Euler equations).
%

Numerical simulations of the model system
exhibit multi-dimensional structures of varying complexity. It was shown that very
regular cellular detonations form for parameter choices near the stability
boundary, and that the farther into the unstable range the parameters are, the more
irregular the patterns become. Characterizing the pattern regularity in the two-dimensional simulations is more challenging than in the corresponding one-dimensional case. Even though transitions from regular to irregular patterns can be observed as parameters of the problem vary, the precise nature of the transitions remains to be further elucidated. In particular, capturing the bifurcation points, and revealing the nature of the apparently irregular solutions (e.g., truly chaotic or not), requires not only extensive long-time simulations, but also more accurate algorithms such as a shock-fitting method, which is known to work well for such problems in one spatial dimension. An interesting problem to explore then is that of the interplay of the one-dimensional and two-dimensional instabilities, especially with regard to the transition from regular to irregular behavior.

The work here is a multi-dimensional extension of a previous 1D analog model
\cite{Faria2014}. Together with \cite{Faria2014}, it demonstrates that much of
the structural and dynamical complexity of gaseous detonations can be captured by
relatively simple models. These models may be of interest in their own right, as they
are very simple, yet produce very complex and rich dynamics. In the context of
reacting flows such as detonations, the present work indicates that rational
asymptotic models of similar complexity may be possible.
Indeed, the asymptotic model in \cite{Faria2014b} (successful in predicting many
properties of multi-dimensional detonations) was derived motivated by ideas such
as the ones here.
Even though the model here is of a semi-\textit{ad hoc} nature,  it
has a number of advantages thanks to its simplicity. In particular, it is easier to solve numerically
and to analyze theoretically (e.g., linear stability analysis) than the full set of reactive Euler equations. Much of the
machinery of detonation theory works for the model, and can be applied
without the heavy algebraic tedium common while working with the reactive Euler
equations. This fact adds a pedagogical value to the analog models, as vehicles
to explain many non-trivial mathematical ideas in detonation theory. In addition,
further advances in understanding may be facilitated by the simplicity of the
analog models.

\appendix

\section{}\label{sec:appendix} 

In this appendix we detail some results that justify the properties stated for the
eigensolutions, $h_j$ and $\theta_j$, introduced in (\ref{eq:eigensol_h}).
Consider a $2\times2$ ODE system of the form
\begin{equation}\label{eq:app01}
 y^\prime = \left(\M - \R\right)\?y\/, \quad 0 < x < \infty\/,
\end{equation}
where $\M$ is a constant matrix with eigenvalues such that $\Re(\mu_1)<\Re(\mu_2)$,
and $\R$ is a matrix valued function of $x$ such that
$\|\R\|_\infty < \infty$ and $\|\R\|_1 < \infty$. Furthermore, we assume that $\M$ and
$\R$ depend on a parameter $\sigma$, analytically in some region $\sigma\in\Omega$
(for notational simplicity, we do not display this dependence below).

\medskip\noindent
Let $\{p_j\}$ be a set of eigenvectors associated with the $\{\mu_j\}$. Then
theorem~8.1 in \cite{Coddington1955} 
shows that a base of solutions $\{y_j\}$ can be selected such that
$y_j \sim e^{\mu_j\?x}\?p_j$ as $x \to \infty$. We will show now that \emph{this base
can also be selected to have analytic dependence for $\sigma\in\Omega$.}

\smallskip\noindent
\emph{Proof.} First of all, $y_1$ is determined uniquely by its asymptotic behavior,
via the Volterra integral equation
\begin{equation}\label{eq:app02}
 z_1 = p_1 + \int_x^\infty e^{(\M-\mu_1)\?(x-s)}\?\R(s)\,z_1(s)\,d\/s = p_1+\calL\?z_1\/,
\end{equation}
where $z_1 = e^{-\mu_1\?x}\?y_1$, and the operator $\calL$ is defined by the formula.
Thus we can write
\begin{equation}\label{eq:app03}
 z_1 = \sum_0^\infty \calL^n\?p_1\/,
\end{equation}
which converges absolutely because
$|\calL^np_1(x)| \leq \frac{1}{n!}\?\left( E\?\|\R\|_1\right)^n\?\|p_1\|\/$,
where $E$ bounds $\left\{e^{(\mu_1-\M)\?t}\/,\, t \geq 0\right\}$. From this it follows
that $z_1$, hence $y_1$, depends analytically on $\sigma$. Then we can take $y_2$
as the solution to (\ref{eq:app01}), with initial data $y_2(0) = z_1(0)^{\perp}$ ---
where $(v_1\/,\,v_2)^{\perp} = (-v_2\/,\,v_1)$. Since $y_2(0)$ depends analytically on
$\sigma$, standard ODE theory \cite{Coddington1955} guarantees that $y_2$ is also
analytic. In addition, because $y_2$ is linearly independent from $y_1$, it must
satisfy $y_2 \sim a\,e^{\mu_2\?x}\?p_2$ for some constant $a \neq 0$.\TheoremEnd


\bibliographystyle{siam}
\bibliography{Bibliography/My-Collection,Bibliography/Xtras}

\begin{thebibliography}{10}

\bibitem{Coddington1955}
{\sc E.~Coddington and N.~Levinson}, {\em {Theory of Ordinary Differential
  Equations}}, Tata McGraw-Hill Education, 1955.

\bibitem{Doering1943}
{\sc W.~D\"{o}ring}, {\em {Uber den detonationvorgang in gasen}}, Annalen der
  Physik, 43(6/7) (1943), pp.~421--428.

\bibitem{Erpenbeck1962}
{\sc J.~J. Erpenbeck}, {\em {Stability of steady-state equilibrium
  detonations}}, Physics of Fluids, 5 (1962), pp.~604--614.

\bibitem{Faria2014thesis}
{\sc L.~M. Faria}, {\em {Qualitative and asymptotic theory of detonations}},
  PhD thesis, KAUST, 2014.

\bibitem{Faria2014}
{\sc L.~M. Faria, A.~R. Kasimov, and R.~R. Rosales}, {\em {Study of a model
  equation in detonation theory}}, SIAM Journal on Applied Mathematics, 74
  (2014), pp.~547--570.

\bibitem{Faria2014b}
\leavevmode\vrule height 2pt depth -1.6pt width 23pt, {\em {Theory of weakly
  nonlinear self sustained detonations}}, arXiv preprint arXiv:1407.8466
  (accepted for publication in the Journal of Fluid Mechanics),  (2015).

\bibitem{Fickett1979}
{\sc W.~Fickett}, {\em {Detonation in miniature}}, American Journal of Physics,
  47 (1979), pp.~1050--1059.

\bibitem{Fickett1985b}
\leavevmode\vrule height 2pt depth -1.6pt width 23pt, {\em {Introduction to
  Detonation Theory}}, University of California Press, 1985.

\bibitem{Fickett2012}
{\sc W.~Fickett and W.~C. Davis}, {\em {Detonation: Theory and Experiment}},
  Courier Dover Publications, 2012.

\bibitem{gamezo1999formation}
{\sc V.~N. Gamezo, D.~Desbordes, and E.~S. Oran}, {\em Formation and evolution
  of two-dimensional cellular detonations}, Combustion and Flame, 116 (1999),
  pp.~154--165.

\bibitem{Hunter2000}
{\sc J.~K. Hunter and M.~Brio}, {\em {Weak shock reflection}}, Journal of Fluid
  Mechanics, 410 (2000), pp.~235--261.

\bibitem{kadomtsev1970stability}
{\sc B.~B. Kadomtsev and V.~I. Petviashvili}, {\em {On the stability of
  solitary waves in weakly dispersing media}}, in Sov. Phys. Dokl., vol.~15,
  1970, pp.~539--541.

\bibitem{Kasimov2013}
{\sc A.~R. Kasimov, L.~M. Faria, and R.~R. Rosales}, {\em {Model for shock wave
  chaos}}, Physical Review Letters, 110 (2013), p.~104104.

\bibitem{Lee2008}
{\sc J.~Lee}, {\em {The Detonation Phenomenon}}, Cambridge University Press,
  2008.

\bibitem{lin1948two}
{\sc C.~C. Lin, E.~Reissner, and H.~S. Tsien}, {\em {On two-dimensional
  non-steady motion of a slender body in a compressible fluid}}, Journal of
  Mathematics and Physics, 27 (1948), p.~220.

\bibitem{Majda1981}
{\sc A.~Majda}, {\em {A qualitative model for dynamic combustion}}, SIAM
  Journal on Applied Mathematics, 41 (1981), pp.~70--93.

\bibitem{may1976simple}
{\sc R.~M. May}, {\em {Simple mathematical models with very complicated
  dynamics}}, Nature, 261 (1976), pp.~459--467.

\bibitem{OranBoris}
{\sc E.~Oran and J.~P. Boris}, {\em {Numerical simulation of reactive flow}},
  Cambridge University Press, Cambridge, UK, 2001.

\bibitem{Oran1998}
{\sc E.~S. Oran, J.~W. Weber, E.~I. Stefaniw, M.~H. Lefebvre, and J.~D.
  Anderson}, {\em {A numerical study of a two-dimensional H2-O2-Ar detonation
  using a detailed chemical reaction model}}, Combustion and Flame, 113 (1998),
  pp.~147--163.

\bibitem{Radulescu2011}
{\sc M.~I. Radulescu and J.~Tang}, {\em {Nonlinear dynamics of self-sustained
  supersonic reaction waves: Fickett's detonation analogue}}, Physical Review
  Letters, 107 (2011), p.~164503.

\bibitem{rosales1989diffraction}
{\sc R.~R. Rosales}, {\em {Diffraction effects in weakly nonlinear detonation
  waves}}, in Nonlinear Hyperbolic Problems, Springer, 1989, pp.~227--239.

\bibitem{sharpe2000two}
{\sc G.~Sharpe and S.~A. E.~G. Falle}, {\em Two-dimensional numerical
  simulations of idealized detonations}, in Proceedings of the Royal Society of
  London A: Mathematical, Physical and Engineering Sciences, vol.~456, The
  Royal Society, 2000, pp.~2081--2100.

\bibitem{tabak1994focusing}
{\sc E.~G. Tabak and R.~R. Rosales}, {\em {Focusing of weak shock waves and the
  von Neumann paradox of oblique shock reflection}}, Physics of Fluids
  (1994-present), 6 (1994), pp.~1874--1892.

\bibitem{vonNeumann1942}
{\sc J.~von Neumann}, {\em {Theory of detonation waves}}, tech. rep., National
  Defense Research Committee Div. B, 1942.

\bibitem{Zab-Khokhlov-1969}
{\sc E.~A. Zabolotskaya and R.~V. Khokhlov}, {\em {Quasi-planes waves in the
  nonlinear acoustics of confined beams}}, Sov. Phys. Acoust., 15 (1969),
  pp.~35--40.

\bibitem{Zeldovich1940}
{\sc Y.~B. Zel'dovich}, {\em {On the theory of propagation of detonation in
  gaseous systems}}, Journal of Experimental and Theoretical Physics, 10
  (1940), pp.~542--569.

\end{thebibliography}

\end{document}